\documentclass[3p]{elsarticle}

\graphicspath{{images/method/}{images/test_elliptic/}{images/test_elasticity/}{images/test_elliptic/revision_misc/}}

\usepackage{amssymb,amsmath,amsthm}
\usepackage{color}
\usepackage{url}
\usepackage{mathrsfs}
\usepackage[ruled]{algorithm2e} 

\newdefinition{rmk}{Remark}

\SetKwProg{Fn}{}{}{end} 

\SetKwComment{Mycomment}{}{}
\SetCommentSty{textit}
\SetKw{And}{and}

\newcommand{\Rset}{\mathbb{R}}

\newcommand{\Nset}{\mathbb{N}}

\usepackage{bm}

\newcommand{\ee}{\mathbf{e}}
\newcommand{\ff}{\mathbf{f}}

\newcommand{\GG}{\mathbf{G}}

\newcommand{\ii}{\mathbf{i}}
\newcommand{\jj}{\mathbf{j}}

\newcommand{\nn}{\mathbf{n}}
\newcommand{\pp}{\mathbf{p}}
\newcommand{\PP}{\mathbf{P}}

\newcommand{\rr}{\mathbf{r}}
\newcommand{\uu}{\mathbf{u}}
\newcommand{\vv}{\mathbf{v}}

\newcommand{\xx}{\mathbf{x}}
\newcommand{\yy}{\mathbf{y}}

\newcommand{\xxi}{\bm{\xi}}

\newcommand{\mcL}{\mathcal{L}}
\newcommand{\mcF}{\mathcal{F}}
\newcommand{\mcB}{\mathcal{B}}

\usepackage[normalem]{ulem}

\newcommand{\La}[1]{#1}
\newcommand{\Lc}[2]{#2}

\begin{document}

\begin{frontmatter}
\title{IGA-based Multi-Index Stochastic Collocation for random PDEs on arbitrary domains}

\author[kaust]{Joakim Beck}
\ead{joakim.beck@kaust.edu.sa}

\author[imati]{Lorenzo Tamellini\corref{cor1}}
\ead{tamellini@imati.cnr.it}

\author[kaust,aachen]{Ra\'{u}l Tempone}
\ead{raul.tempone@kaust.edu.sa}

\cortext[cor1]{Corresponding author}
\address[kaust]{King Abdullah University of Science and Technology (KAUST), Computer, Electrical and Mathematical Science and Engineering Division (CEMSE)
  Thuwal 23955-6900, Saudi Arabia}
\address[imati]{Consiglio Nazionale delle Ricerche - Istituto di Matematica Applicata e Tecnologie Informatiche
  ``E. Magenes'' (CNR-IMATI), Via Ferrata 1, 27100, Pavia, Italy}
\address[aachen]{Alexander von Humboldt Professor in Mathematics of Uncertainty Quantification, RWTH Aachen University, 52062 Aachen, Germany}
\begin{abstract}
This paper proposes an extension of the Multi-Index Stochastic Collocation (MISC) method for forward uncertainty
quantification (UQ) problems in computational domains of shape other than a square or cube, by exploiting
isogeometric analysis (IGA) techniques. Introducing IGA solvers to the MISC algorithm
is very natural since they are tensor-based PDE solvers, which are precisely what is required by the MISC machinery.
Moreover, the combination-technique formulation of MISC allows the straight-forward reuse of
existing implementations of IGA solvers. We present numerical results to showcase the effectiveness of the proposed approach.
  
  \medskip
  \textbf{Highlights}
  \begin{itemize}
  \item Isogeometric solvers used in a MISC framework for forward UQ problems.
  \item The combination-technique formulation of the method allows straight-forward reuse of legacy IGA solvers.
  \end{itemize}

\end{abstract}

\begin{keyword}
  Isogeometric analysis \sep Uncertainty Quantification \sep Sparse Grids \sep Stochastic Collocation methods \sep
  multilevel methods \sep combination-technique
\end{keyword}

\end{frontmatter}




\section{Introduction}

Uncertainty quantification (UQ) has received a considerable amount of attention in recent years and is by now
considered an essential tool in the domain of computational science and engineering \cite{ghanem:UQbook,sullivan:UQbook,smith:UQbook}.
However, performing UQ analyses remains a significant computational challenge, since these analyses typically require
solving a computational model repeatedly for different values of the uncertain variables in the model.
In this paper, we consider in particular models described by elliptic PDEs whose solution is denoted by $\uu$.
Two general ``meta-strategies'' (complementary to each other) are by now recognized in the UQ community as key
to reducing the computational cost and making UQ analyses feasible:
a) dimension-adaptivity, i.e., investing the majority of the computational cost in approximating
the dependence of $\uu$ on the random variables whose variability has the largest impact on $\uu$ itself; and
b) a multi-level approach, in which most of the variability of $\uu$ is explored by using ``low-fidelity''
approximations of the PDE (e.g., coarse meshes and simplified-physics models), and resorting
to ``high-fidelity'' approximations only sparingly.
Dimension-adaptivity was the first ``meta-strategy'' to be introduced in the UQ community and has been
extensively discussed in the literature, e.g. in
\cite{nobile.tempone.eal:aniso,nobile.eal:optimal-sparse-grids,cohen.devore.schwab:nterm2,chkifa:adaptive-taylor},
while the multi-level approach was developed more recently; see, e.g.,
\cite{scheichl.giles:MLMC,hajiali.eal:MultiIndexMC,peherstorfer:MFsurvey}.

The Multi-Index Stochastic Collocation (MISC) method for random PDEs was first introduced in \cite{hajiali.eal:MISC1,hajiali.eal:MISC2}
and represents an attempt to combine both strategies. MISC uses an ad-hoc algorithm
to simultaneously choose both the best sequence of computational meshes and the random variables
whose impact should be more carefully approximated. It is closely related to the 
sparse-grids technique introduced to solve high-dimensional PDEs
\cite{Bungartz.Griebel.Roschke.ea:pointwise.conv,b.griebel:acta,Griebel.schneider.zenger:combination,Hegland:combination}
and to compute high-dimensional integrals/interpolants \cite{smolyak:quadrature,barthelmann.novak.ritter:high},
and indeed it can be seen as a combination of the two methodologies.
The main idea behind the MISC algorithm is to write the approximation operator
as a linear combination of many less-refined approximations, in the spirit of a Richardson extrapolation \cite{quarteroni.sacco.eal:numerical}.
In particular, a ``profit'' is assigned to each possible component, and eventually only those components with the largest profit
are included in the computation.
Another possible approach would be to balance the errors of the physical and stochastic discretizations,
as proposed in \cite{teckentrup.etal:MLSC,van-wyk:MLSC,hps13}.

Crucially, MISC needs both a physical solver and a sampler in the stochastic domain with a tensor structure.
The tensor construction in the stochastic domain can be obtained by e.g. tensorizing standard Lagrangian (interpolant) polynomials,
which would also be the starting point of the classical sparse-grids collocation method for UQ
\cite{babuska.nobile.eal:stochastic2,xiu.hesthaven:high}.
When it comes to the physical solver, previous works on MISC \cite{hajiali.eal:MISC1,hajiali.eal:MISC2}
only considered square domains over structured grids, which induce a tensor structure in the solver;
possible strategies to apply MISC to non-square domains, on which having a tensorized solver might be non-trivial,
were only briefly mentioned in \cite{hajiali.eal:MISC2}.

In this paper, we extend \cite{hajiali.eal:MISC1,hajiali.eal:MISC2} and consider MISC on arbitrarily shaped domains,
by employing isogeometric analysis (IGA) \cite{Hughes:2005,IGA-book}.
IGA is an alternative to the standard finite element analysis that uses
the basis functions employed by CAD softwares to represent the computational domain
(typically, B-splines or non uniform rational B-splines (NURBS)) as basis for the approximation of the PDE solution as well.
IGA has several interesting features and has therefore received a growing
interest from researchers and practitioners in computational science.
For example, IGA can work with exact domain representation and the meshing process is simplified in some situations;
B-splines/NURBS of arbitrary polynomial degree and regularity can be generated in a very easy way 
and show, in certain cases, superior error vs. degrees-of-freedom ratio with respect to standard finite element
bases. We refer the reader to \cite{acta-IGA} for an in-depth discussion on IGA. Crucially, multivariate B-splines and NURBS
are built by tensorization of their univariate counterparts, which makes IGA solvers
particularly suitable for use with MISC.
\La{We remark that not many works focusing on UQ are available in the IGA literature,
  and to the best of our knowledge none of them considers multi-level/multi-index strategies; here we mention
  \cite{Benzaken20171215,manzoni.heltai:RB-IGA,WILHELM2016,corno_UQIGA,LI_UQIGA,HIEN2017,RAHMAN2018}.}

It is important to remark that in the end the MISC algorithm itself simply prescribes
to solve a number of standard uncoupled PDEs, each of them corresponding to a different realization
of the random parameters, on physical meshes with different resolutions (possibly
anisotropic, i.e., more refined along some physical directions). Therefore,
any available IGA software can be readily re-used. 
We also point out that in the case where there are no random variables,
this procedure corresponds to the sparse IGA method for solving PDEs
discussed in \cite{beck.eal:sparse-IGA}.
Finally, we mention that while IGA is a convenient choice to extend MISC to problems on
non-square domains, it is not the only possible choice. Other choices, such as finite differences, finite volumes, $\mathbb{Q}_k$ finite elements,
might be envisaged; the comparison of these methods is outside the scope of this work.

\bigskip

The rest of this paper is organized as follows. We introduce the
general UQ framework for elliptic PDEs with random coefficients in Section \ref{section:problem_def},
and discuss IGA solvers in Section \ref{section:IGA}. We present the MISC algorithm
in Section \ref{section:MISC} and showcase the results obtained with MISC on
some numerical examples in Section \ref{section:results}.
Finally, we draw conclusions in Section \ref{section:conclusions}.

\bigskip

Throughout the manuscript, we make extensive use of multi-indices, i.e.,
vectors with integer components. To this end, we recall some useful definitions and notations:
\begin{itemize}
\item Given $\ii,\jj \in \Nset^K$, $\ii \leq \jj$ means that $i_k \leq j_k$ for $k=1,\ldots,K$;
\item $\bm{e}_i$ is the $i$-th canonical multi-index, i.e.,   $(\bm{e}_i)_{k}=1$ if $i=k$, and zero otherwise;
\item $\bm{1}$ is a multi-index whose components are all equal to 1;
\item Given a function $f:\Rset\rightarrow\Rset$ and a multi-index $\ii$,
  $\mathbf{f}(\ii)$ denotes the multi-index $[f(i_1),f(i_2),\ldots,f(i_K)]$;
\item A multi-index set $\Lambda \subset \Nset^K$ is said to be downward closed if
  \begin{equation}\label{eq:downward_closed_set}
  \forall \ii \in \Lambda, \quad \ii - \ee_j \in \Lambda \text{ for every } j=1,\ldots,K \text{ such that } \ii_j > 1.  
  \end{equation}
  
\item The margin of a multi-index set $\Lambda$, $\text{Mar}(\Lambda)$, is the set of multi-indices that can be reached
  ``within one step'' from $\Lambda$,
  \begin{equation}\label{eq:margin}
    \text{Mar}(\Lambda) = \{\ii \in \Nset^K \text{ s.t. } \ii = \jj + \ee_k \text{ for some } \jj \in \Lambda \text{ and some } k \in \{1,\ldots,K\} \};    
  \end{equation}
\item The reduced margin of  a multi-index set $\Lambda$, $\text{Red}(\Lambda)$, is the subset of $\text{Mar}(\Lambda)$ composed
  by indices from which ``every backward step'' will take inside $\Lambda$,
  \begin{equation}\label{eq:reduced_margin}
    \text{Red}(\Lambda) = \{\ii \in \Nset^K \text{ s.t. } \ii - \ee_k \in \Lambda \text{ for every } k \in \{1,\ldots,K\} \text{s.t. } i_k > 1 \}.    
  \end{equation}
\end{itemize}

\section{Problem definition}\label{section:problem_def}

Let $\mcB$ be a compact domain in $\Rset^d$, $d=2,3$, that represents the ``physical domain'' of the problem.
In addition, let $\yy =[y_1,y_2,\ldots,y_N]$ be an $N$-dimensional random vector, whose components are mutually independent
random variables with support $\Gamma_n \subset \Rset$ and probability density function $\rho_n(y_n)$. 
Thus, $\yy \in \Gamma = \Gamma_1 \times \Gamma_2 \cdots \Gamma_N$, and $\Gamma$ represents the ``stochastic domain'' of the problem;
since $y_n$ are mutually independent, $\rho(\yy) = \prod_{n=1}^N \rho_n(y_n)$ is a probability density function on $\Gamma$.
Throughout this work, we often refer to $y_n$ as ``stochastic directions'', which is a short-hand for ``directions of the stochastic domain''.

We consider the following problem: Find $\uu : \mcB \times \Gamma \to \Rset^m$ such that for $\rho$-almost every $\yy \in \Gamma$,
\begin{equation}\label{eq:strong-form}
\begin{cases}
\mcL(\uu;\xx,\yy)=\mcF(\xx) & \xx \in \mcB, \\
\uu(\xx,\yy)=0 & \xx \in \partial\mcB,
\end{cases}
\end{equation}
where $\mcL$ is a differential operator and $\mcF$ is an operator on $\xx$.
In particular, in the numerical results section we consider a linear scalar elliptic equation and
a linear elasticity equation. We assume well-posedness of the problem in some Hilbert space $V$
for $\rho$-almost every $\yy \in \Gamma$  (specific choices of $V$ are detailed for each example in Section \ref{section:results});
observe that $\uu$ can also be seen as an $N$-variate Hilbert-space-valued function $\uu(\yy): \Gamma \to V$,
and in particular it is convenient to introduce the Bochner space of finite-variance Hilbert-space-valued functions, 
$L_{\rho}^2(\Gamma;V)=\{\uu: \Gamma \to V \ \text{strongly measurable such that} \int_{\Gamma} \Vert \uu(,\cdot,\yy) \Vert^2_{V} \rho(\yy)d\yy < \infty \}$,
to which $\uu$ is assumed to belong. 

The random variables $\yy$ model the uncertainties in the system, i.e., account for the fact that coefficients, forcing terms,
boundary/initial conditions, and domain shape are often ``imperfectly'' known due to measurement errors,
lack of data, or intrinsic variability (e.g., when they describe phenomena like wind, rain, earthquakes). The goal of a forward UQ
analysis is therefore to assess how much the variability of such random objects affects the quantities of interest
of the computation, which could be either the solution $\uu$ or a functional thereof. In particular, in this work
we assume that some functional of the solution $\uu$, $\Phi: V \to \Rset$,
e.g. $\Phi(\vv) = \int_{\mcB} \vv(\xx) d\xx$ or $\Phi(\vv) = \vv(\xx_0)$  is given, and we aim to estimate its expected value,
i.e., we want to compute
\[
\mathbb{E}[\Phi(\uu(\xx,\yy))]=\int_{\Gamma}\Phi(\uu(\xx,\yy))\rho(\yy)d\yy.  
\]
In this work, we numerically analyze the performance of MISC for this task.
More specifically, we use MISC to select the physical and stochastic discretizations parameters,
but other discretization parameters (e.g., time-steps, number of particles,
solver tolerances) could be added to the set of parameters governed by MISC.
The novelty of the present work consists in showing how to extend the MISC methodology to more general physical domain shapes,
by replacing the multi-linear finite elements solver adopted in previous works with the IGA solvers that
we present in the next section.

\section{Isogeometric solvers}\label{section:IGA}

In this section, we briefly present the fundamentals of isogeometric analysis (IGA)
and refer the reader to \cite{Hughes:2005,IGA-book,acta-IGA} for a more thorough discussion.
The first ingredient of IGA is a \Lc{B-splines/NURBS}{} representation \La{over a B-splines/NURBS basis}
of the computational domain $\mcB$, which is usually provided by a CAD software.
An isogeometric solver typically then uses the same set of basis functions to compute an approximation of the solution of the PDE;
\La{this is the approach that we consider in this work.}
However, \Lc{that}{} this is not strictly needed, and two different sets of B-splines/NURBS functions could be used instead: one to
approximate the geometry and another to approximate the solution.

The B-splines/NURBS representation of $\mcB \subset \Rset^d$ consists of a transformation from a reference
domain $\widehat{\mcB}$ (typically, a square or cube)
to the physical domain $\mcB$, written as a linear combination of B-splines/NURBS functions
and so-called control points $\PP_\ii \in \Rset^d$. Following the IGA/CAD literature, we refer in this
work to $\widehat{\mcB}$ as the ``parametric domain''. 

\begin{figure}
  \centering
  \includegraphics[width=0.32\linewidth]{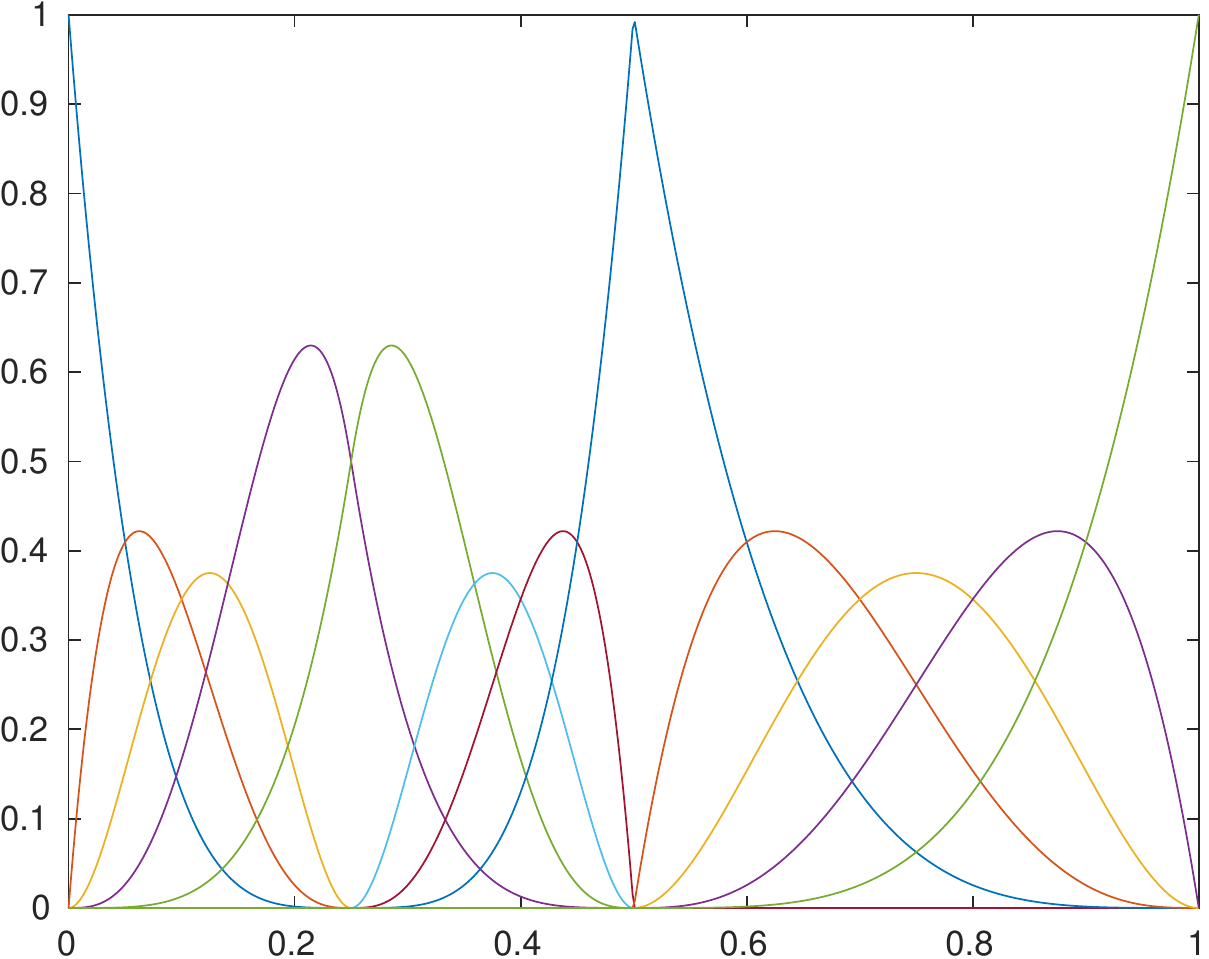}
  \includegraphics[width=0.32\linewidth]{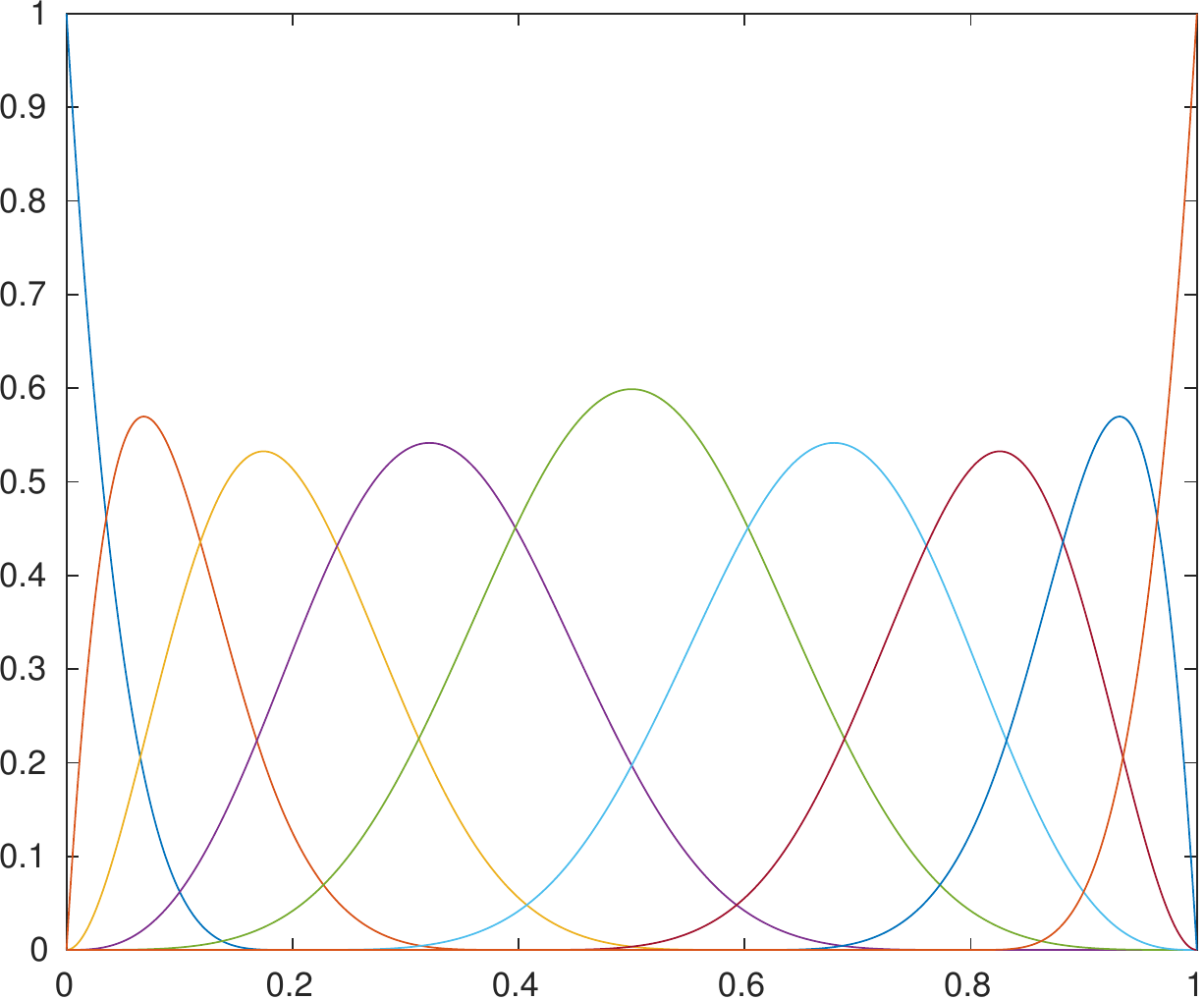}
  \includegraphics[width=0.32\linewidth]{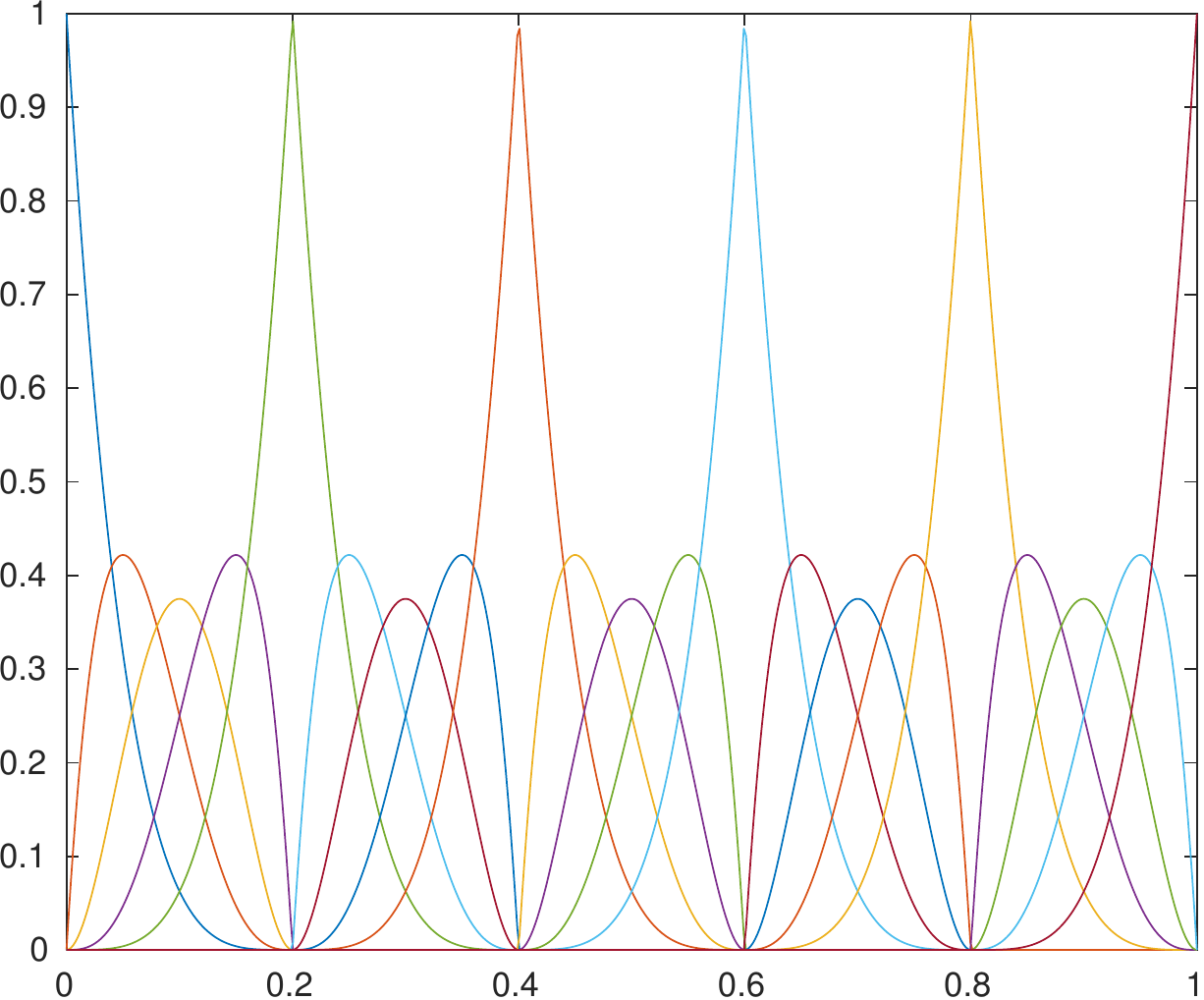}
  \caption{B-splines basis built on different knot vectors; each knot vector is ``open'', i.e.,
    the first and last nodes are repeated $p+1$ times, with $p=4$.
    Left: $\Xi= [0,0,0,0,0,0.25,0.25,0.25,0.5,0.5,0.5,0.5,1,1,1,1,1]$.
    Note the different continuity at $\xi=0.25$ and $\xi=0.5$ due to different multiplicity:
    $\xi$ is repeated 3 times, thus the B-splines are $C^1$ in that knot;
    $\xi$ is repeated 4 times, thus the B-splines are $C^0$ in that knot.
    Center: $\Xi=[0,0,0,0,0.2,0.4,0.6,0.8,1,1,1,1,1]$. In this case, the internal knots are repeated only once,
    resulting in maximally smooth B-splines, i.e., $C^3$ at the internal knots.
    Right: $\Xi=[0,0,0,0,0,0.2,0.2,0.2,0.2,0.4,0.4,0.4,0.4,0.6,0.6,0.6,0.6,0.8,0.8,0.8,0.8,1,1,1,1,1]$.
    Internal knots are repeated $p$ times, resulting in $C^0$ splines (i.e., Lagrangian finite elements of degree $p$).
    As expected, the higher the continuity of the basis functions, the smaller the cardinality of the basis.}
  \label{fig:Bsplines_constr}
\end{figure}

The B-splines/NURBS used in such representation are built by tensorization, 
and therefore we begin our presentation by considering the univariate case.
We introduce a reference interval, $\hat{I}=[0,1]$, and a knot vector over $\hat{I}$,
i.e, a non-decreasing vector $\Xi=[\xi_{1},\xi_{2}...,\xi_{n+p+1}]$, where $n,p \in \Nset$
and $\xi_1, \xi_{n+p+1}$ coincide with the extrema of $\hat{I}$,
that will be used to define a set of $n$ B-splines polynomials of degree $p$ over $\hat{I}$;
observe that $\Xi$ can have repeated entries, for reasons that will be clear later on (see Figure \ref{fig:Bsplines_constr}).
Each $\xi_i$ is called ``knot'', and an interval $(\xi_i, \xi_{i+1})$ of non-zero length is an
``element''; $N_{el}$ is the number of elements. The elements are not required to have the same length,
but if they do, we call that length the mesh-size, which is denoted by $h$.
We define the non-decreasing vector $Z= [\zeta_1, \ldots, \zeta_{N_{el}+1}]$ as the vector of knots of $\Xi$ without
repetitions, and $m_i$ is the multiplicity of $\zeta_i$ in $\Xi$, such that $\sum_{i=1}^{N_{el}} m_i = n+p+1$. 
A knot vector is said to be ``open'' if its first and last knots have a multiplicity of $p+1$.

We can now define the B-splines polynomials of degree $p$ on $\hat{I}$ by means of the Cox-De Boor recursive formula \cite{acta-IGA}:
we start with piecewise constant ($\tilde{p}=0$),
\[
\widehat{S}_{i,0}(\xi)= 
\begin{cases}
  1 &  \xi_{i}\leq \xi<\xi_{i+1} \\ 
  0 &  \textrm{otherwise,} \\
\end{cases} \qquad \qquad  \mbox{for } i = 1,\ldots,n+p
\]
and then, for $\tilde{p} = 1,\ldots,p$, we have the recursive step 
\[
  \widehat{S}_{i,\tilde{p}}(\xi)=
  \begin{cases}
    \dfrac{\xi-\xi_{i}}{\xi_{i+\tilde{p}}-\xi_{i}}\widehat{S}_{i,\tilde{p}-1}(\xi)+\dfrac{\xi_{i+\tilde{p}+1}-\xi}{\xi_{i+\tilde{p}+1}-\xi_{i+1}}\widehat{S}_{i+1,\tilde{p}-1}(\xi),
       & \xi_{i}\leq \xi<\xi_{i+\tilde{p}+1} \\
       0, & \textrm{otherwise}
     \end{cases}
     \quad \mbox{ for } i = 1,\ldots,n+p-\tilde{p},
\]  
with the understanding that $0/0=0$; note that if the knot vector $\Xi$ is open, the corresponding basis
will be interpolatory in the first and last knots. The B-splines are polynomials of degree $p$
and continuity $C^{p-m_i}$ at $\zeta_i$, which means that the regularity of
the B-splines can be reduced by repeating multiple times the same entry of the knot vector. In particular,
repeating a knot $p$ times will result in a basis with $C^0$ regularity in that knot, which means that
the basis will also be interpolatory at that knot; see Figure \ref{fig:Bsplines_constr}-left
for an example. The generated B-splines are linearly independent, and we refer to their span as the ``space of splines'', i.e., $W_p(\Xi,\hat{I})$:
\[
W_p(\Xi,\hat{I}) = \text{span}\left\{ \widehat{S}_{i,p}\,,\,\, i=1,\ldots,n \right\}.
\]
Often, all the internal knots are repeated the same number of times to obtain a B-splines basis with
continuity $r$ at each point $\zeta_i$, $0\leq r \leq p-1$. If that is the case, we add a
superscript $r$ to the notation, i.e., $S_{i,p}^r$ or $W_p^r(\Xi,\hat{I})$. Clearly, for a fixed polynomial degree $p$,
the number of basis functions $n$ decreases as $r$ increases; see Figure \ref{fig:Bsplines_constr}-center and right
for examples of $C^{p-1}_p$ and $C^0_p$ B-splines.

Next, we define B-splines on $d$-dimensional domains for the case $d=2$, with the understanding that the extension to the case $d=3$ is trivial.
The fact that multivariate B-splines are defined by tensorization is crucial for the development
of the MISC methodology, as will become clear in Section \ref{section:MISC}.
We define the parametric domain $\widehat{\mcB} = \hat{I} \times \hat{I}$ and consider two
open knot vectors $\Xi_1, \Xi_2$ with $n_1 + p_1 + 1$ and $n_2+p_2+1$ knots, respectively;
the corresponding knots without repetitions are denoted by $Z_1,Z_2$.
We introduce the tensor products $\mathbf{\Xi} = \Xi_1 \otimes \Xi_2$ and $\mathbf{Z}=Z_1 \otimes Z_2$;
$\mathbf{Z}$ generates a Cartesian mesh over $\widehat{\mcB}$ 
composed by $N_{el,1} \times N_{el,2}$ rectangular elements.
According to the CAD/IGA literature we refer throughout this work to the two directions $\xi_1$ and $\xi_2$
as ``parametric directions''.
A basis for the space of bi-variate splines is then obtained by
taking tensor products of the univariate B-splines, 
\[
W^\rr_\pp(\mathbf{\Xi},\hat{\mcB}) = \text{span}\{\widehat{S}_{\ii,\pp}, \ii \leq \nn \},
\]
where $\ii=[i_1,i_2]$, $\pp=[p_1,p_2]$, $\nn=[n_1,n_2]$, $\rr=[r_1,r_2]$,
and $\widehat{S}_{\ii,\pp}(\xi_1,\xi_2) = \widehat{S}_{i_1,p_1}(\xi_1) \widehat{S}_{i_2,p_2}(\xi_2)$. 

\begin{figure}[thb]
  \centering
  \includegraphics[width=0.32\linewidth]{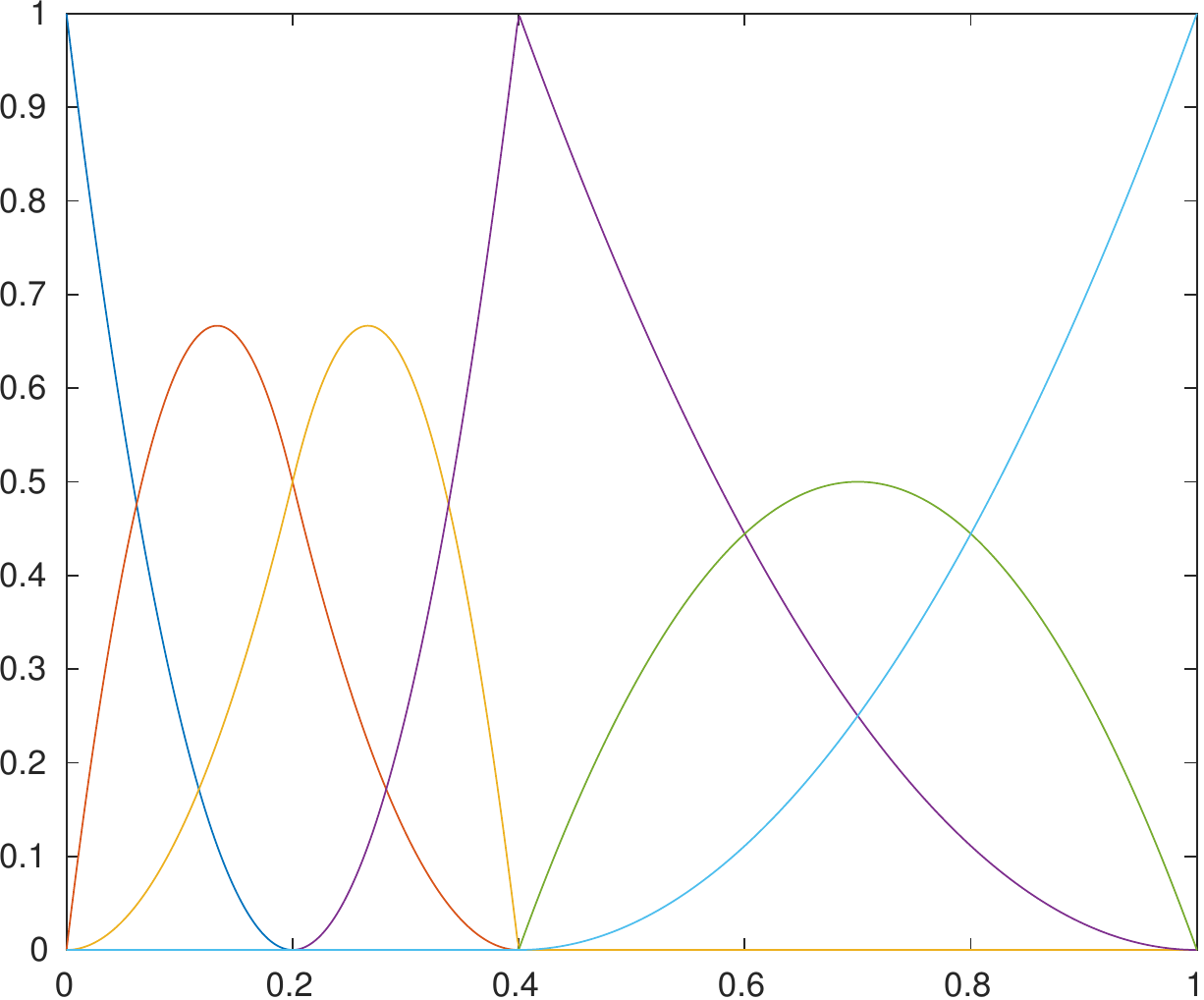}
  \includegraphics[width=0.32\linewidth]{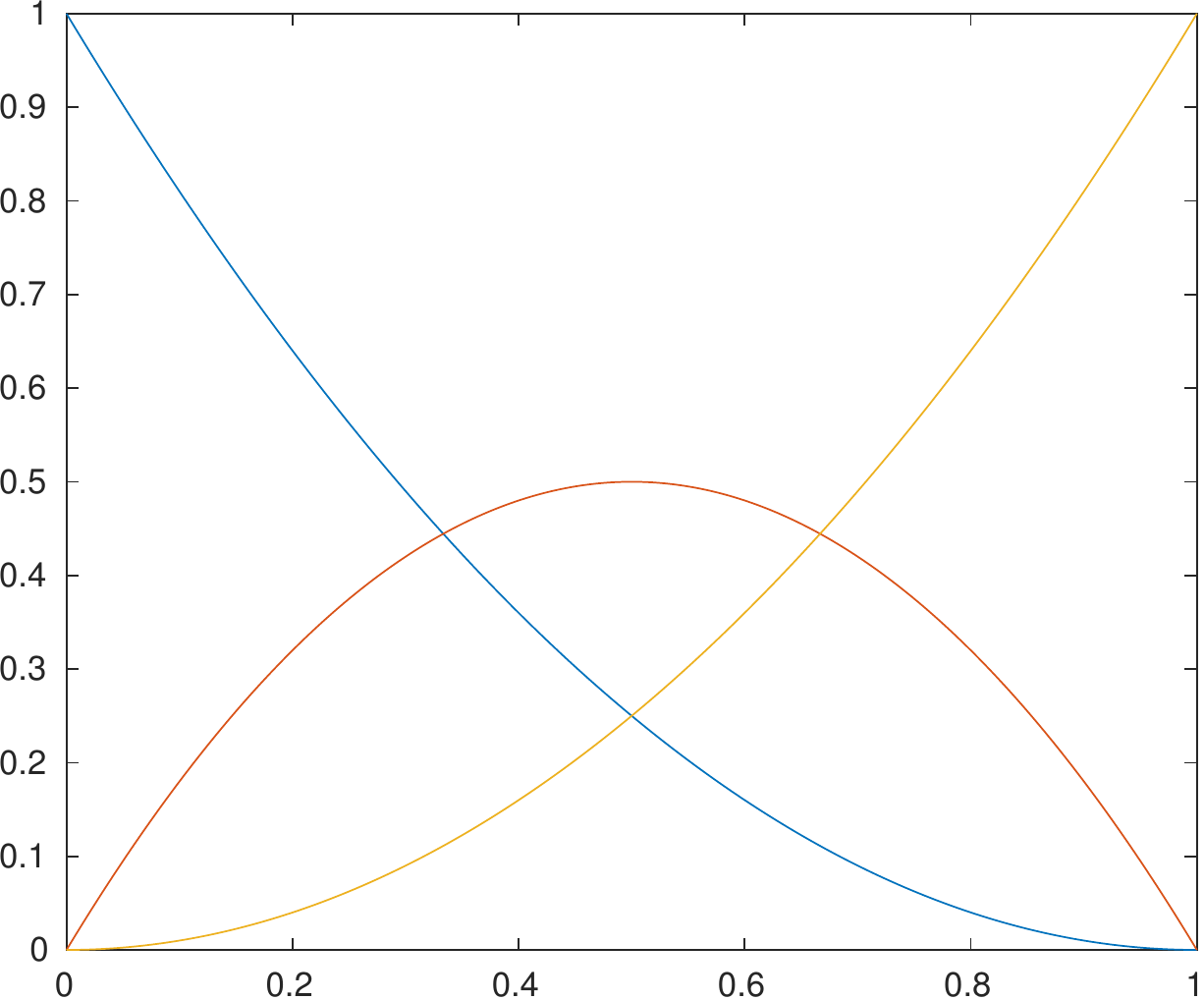}
  \includegraphics[width=0.32\linewidth]{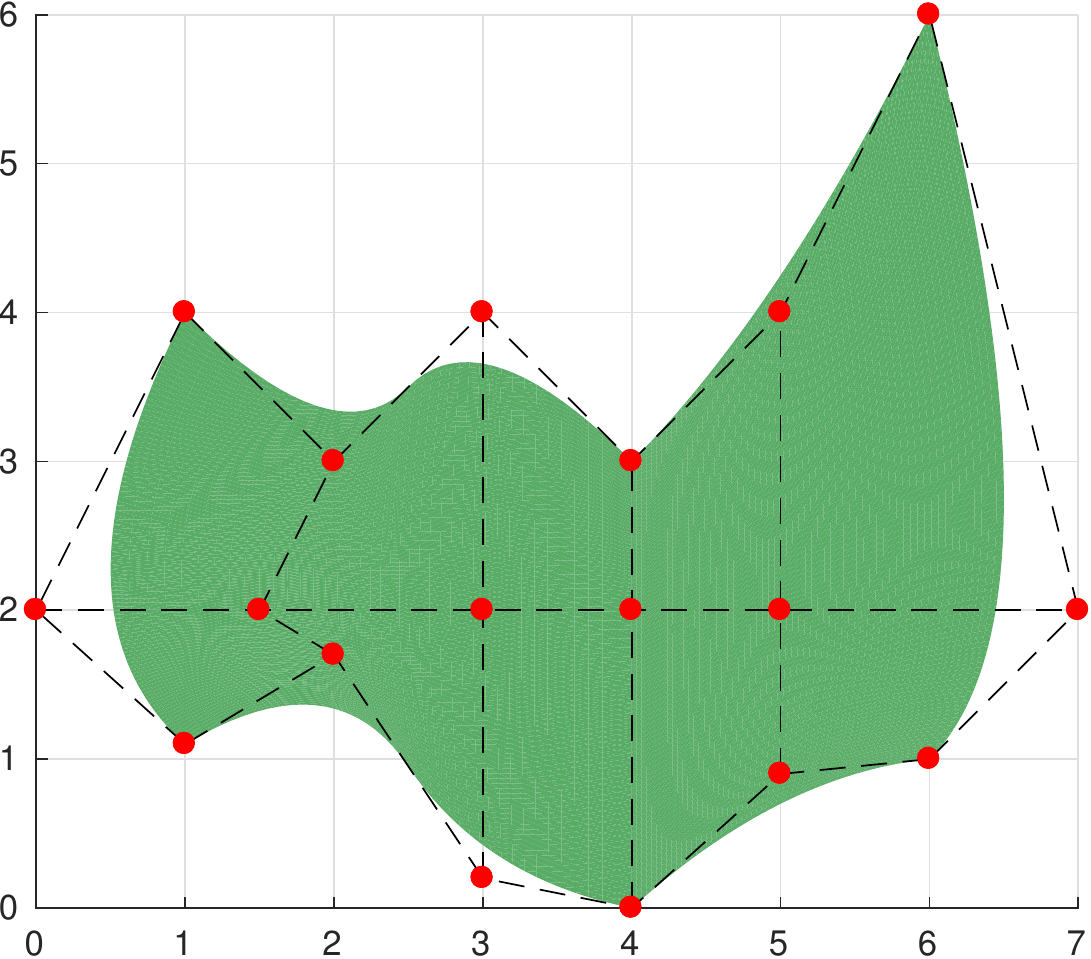}
  \caption{Surfaces are built by linear combinations of control points and tensorized B-splines.
    Left: B-spline basis along the first reference interval;
    center: B-spline basis along the second reference interval. The tensor product of these two bases forms
    the basis (with cardinality 18) of the linear combination used to represent the surface in the right plot.
    The corresponding 18 control points \La{$\PP_{\ii}$} are marked in red. The dotted lines connecting the control points represent
    the so-called ``control net''.}
  \label{fig:planar_surf}
\end{figure}

We are now in the position to introduce the B-splines representation of the computational domain $\mcB$ using
a linear combination of B-splines with control points $\PP_\ii \in \Rset^2$ and $\ii \leq \nn$ 
(see also Figure \ref{fig:planar_surf}):
\[
\xx \in \mcB \Leftrightarrow \xx = \GG(\xxi) = \sum_{\ii \leq \nn} \PP_{\ii} \widehat{S}_{\ii,\pp} (\xxi) \text{ for some }\xxi \in \widehat{\mcB}.  
\]
In the CAD literature, the function $\GG:\widehat{\Omega} \rightarrow \Omega$ is often called a ``parameterization'' of the geometry $\mcB$;
\La{the control points $\PP_\ii$ are chosen by the CAD designer during the design phase, see e.g. \cite{farin2001curves,cohen2001geometric}}.
Observe that the control points need not belong to $\mcB$: this is the case only if the corresponding
basis function is $C^0$ continuous; 
see again Figure \ref{fig:planar_surf}. 
With a slight abuse of wording, in the following we sometimes talk of ``physical directions''
instead of ``parametric directions''; we also use ``physical directions'' as a shorthand for the longer
``curvilinear coordinates induced by mapping the parametric directions over the physical space''.
However, this talk of ``physical directions'' will be useful when juxtaposed with the ``stochastic directions'' $y_1,\ldots,y_N$.

\begin{figure}[thb]
  \centering
  \includegraphics[width=0.42\linewidth]{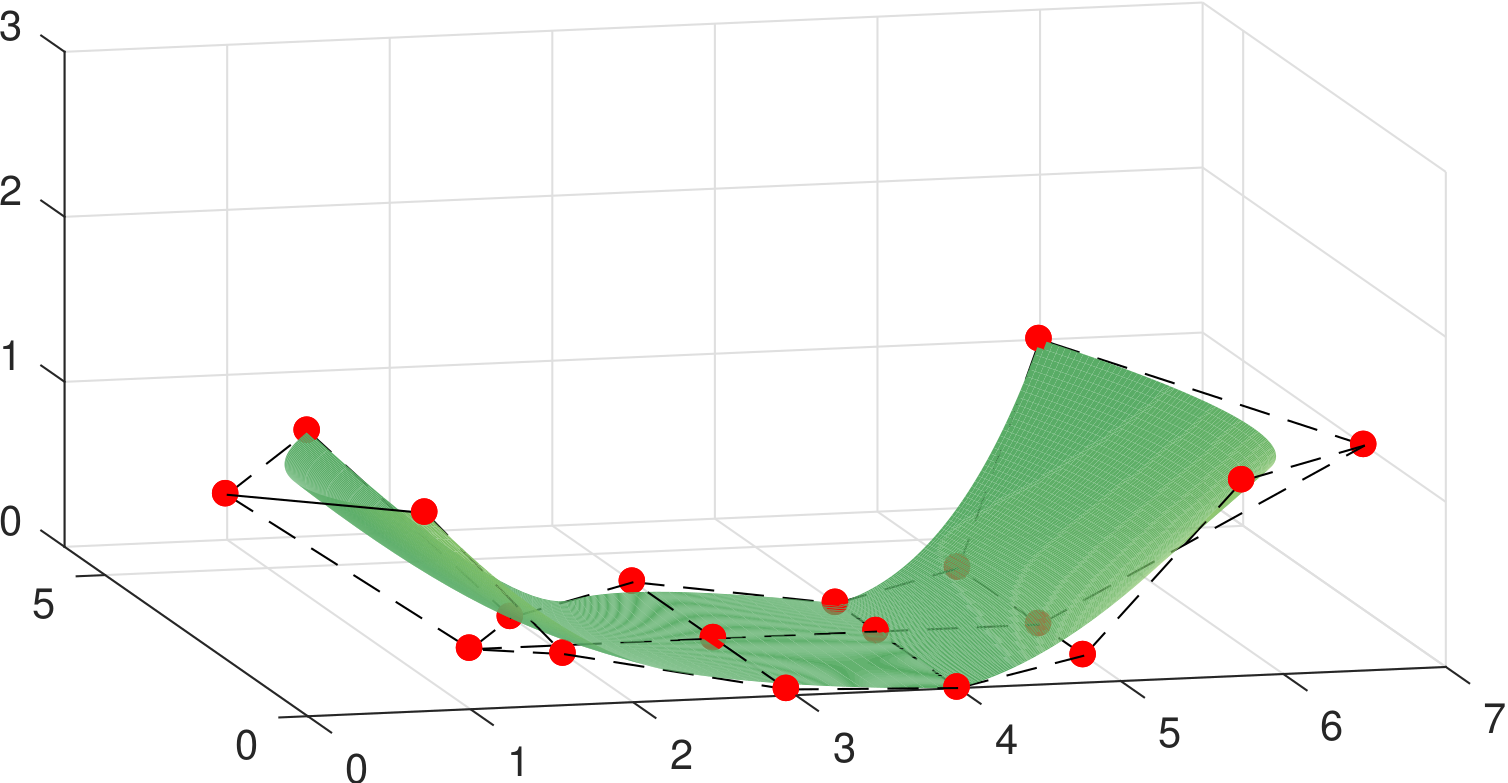}
  \includegraphics[width=0.32\linewidth]{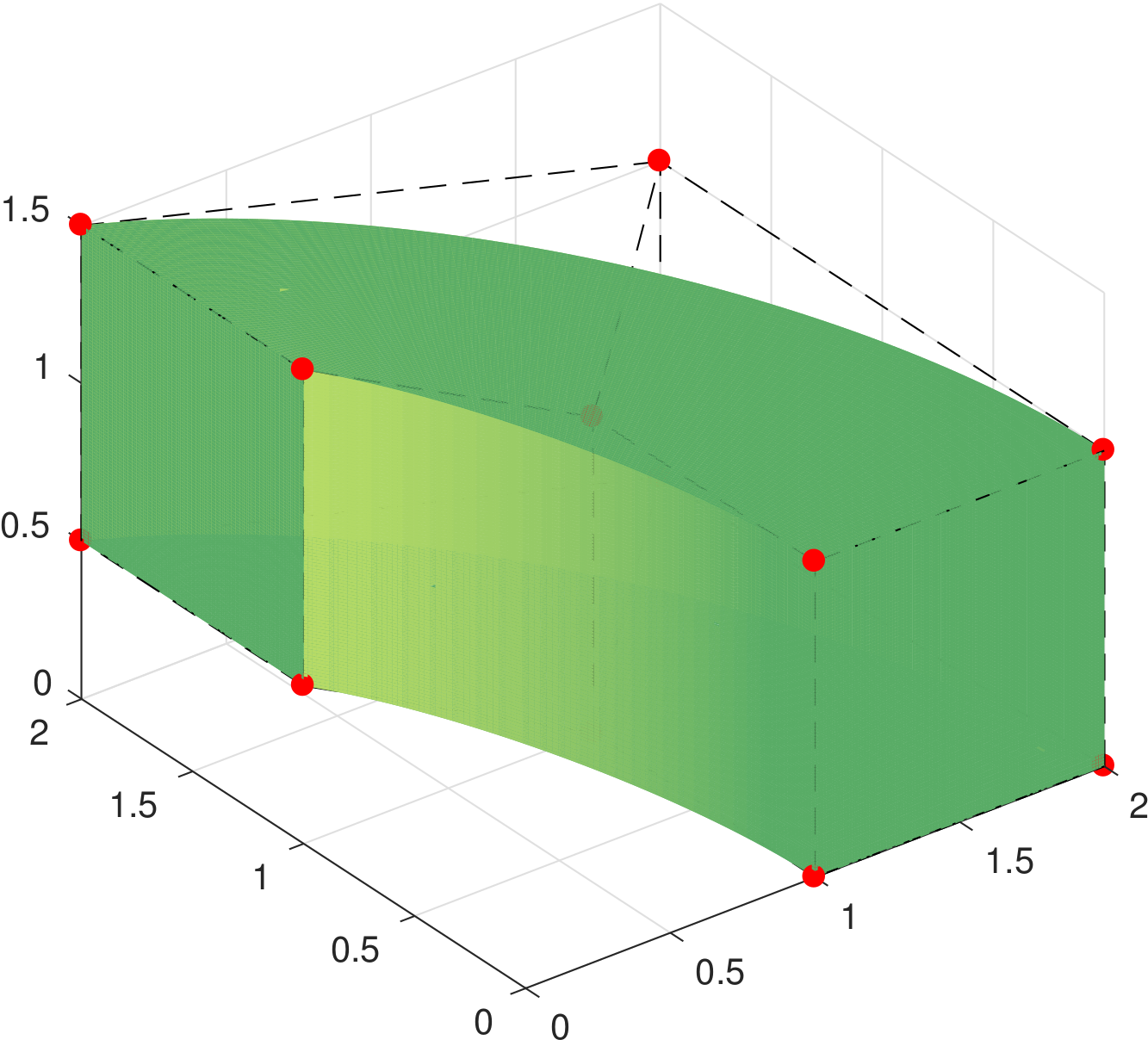}
  \caption{Left: A non-planar surface, obtained from the example in Figure \ref{fig:planar_surf} by choosing control points in $\Rset^3$ instead of $\Rset^2$.
    Right: A volume is represented by taking tensor products of three B-splines bases.}\label{fig:nonplanar_and_volumes}
\end{figure}

Non-planar surfaces can be also generated in the same way, by choosing $\PP_\ii \in \Rset^3$ instead of $\PP_\ii \in \Rset^2$;
see Figure \ref{fig:nonplanar_and_volumes}-left.
A further round of tensorizations, again with $\PP_\ii \in \Rset^3$, would allow us to represent volumetric computational
domains $\mcB$; see Figure \ref{fig:nonplanar_and_volumes}-right. 
Incidentally, note that many geometries of practical interest, and in particular all conic sections other than parabulae,
cannot be represented exactly by B-splines. To this end, the so-called 
\emph{non-uniform rational B-splines} (NURBS) have been introduced, see again \cite{Hughes:2005,IGA-book,acta-IGA}
for details. As the name suggests, NURBS are ratios of B-splines and retain most of the properties of B-splines.

In IGA, the B-splines basis is also used to approximate the solution to \eqref{eq:strong-form}.
We therefore introduce the B-splines space on the physical domain $\mcB$ as
\[
  W_{\pp}(\mathbf{\Xi},\mcB) = \text{span} \{S_{\ii,\pp} = \widehat{S}_{\ii,\pp} \circ \GG^{-1}, \ii \leq \nn  \},
\]
and then approximate the solution to \eqref{eq:strong-form} for a fixed value of the random variables $\yy$  as
\begin{equation}\label{eq:u-IGA-approx}
  \uu(\xx,\yy) \approx \uu_\nn(\xx,\yy)  = \sum_{\ii \leq \nn} c_{\ii} S_{\ii,\pp}(\xx),
\end{equation}
where the coefficients $c_\ii$ depend on $\yy$ (for simplicity, we do not show this dependence in the equation).
The coefficients $c_\ii$ can be computed either by following a standard Galerkin approach or by collocation methods,
i.e., inserting \eqref{eq:u-IGA-approx} in the strong form of the problem, \eqref{eq:strong-form},
and enforcing that the equation be satisfied in a set of collocation points, which results 
in a system of linear equations to be solved \cite{montardini:collocation,gomez2016variational,anitescu2015isogeometric,casquero2016isogeometric}.
In this work, we employ a standard Galerkin method; we remark again that,
as will be clearer later, pre-existing solvers can be readily reused in the context of the MISC method.

As expected, the approximation of $\uu$ in \eqref{eq:u-IGA-approx}, $u_{\nn}$,
converges to $\uu$ as the cardinality of $W_\pp(\mathbf{\Xi},\mcB)$ increases,
and the convergence results with respect to both $h$ and $p$ are analogous to the standard finite elements results.
\La{We conclude this introduction to IGA pointing out that IGA is still in its infancy, and as such quite many issues
  remain open and limit to some extent the applicability of the IGA approach so far. Here we mention only two such aspects,
  which both concern the interplay between geometry representation and PDE solutions.
  First, it is important to point out that most industrial CAD software represent 3D objects by their boundaries only,
  so that a preliminar non-trivial conversion to a 3D B-Spline/NURBS model is needed \cite{DOKKEN2018,MASSARWI2016};
  moreover, such software often complex geometries are represented by
  boolean operation between elementary shapes by so-called ``trimming operations'', which partially destroy the tensor
  structure and require ad-hoc treatment \cite{KUDELA2016406,Marussig2018}.}

In the following, we fix the polynomial degrees along each parametric direction, i.e., $p_1=\ldots=p_d=p$,
and consider $h$-refinements only, i.e., we increase the cardinality of the basis by adding knots
to the initial knot vectors $\Xi_1$ and $\Xi_2$, though it is also possible to devise a $p$-MISC technique (or a combined $h-p$ version).
In particular, for each parametric direction $\xi_1, i=1,\ldots,d$, we consider a sequence of knot vectors $\Xi_i$
indexed by an integer $\alpha_i \geq 1$, whose number of elements doubles whenever $\alpha_i$ increases by 1: $N_{el,i} \sim 2^{\alpha_i}$.
We let $\bm{\alpha} = [\alpha_1,\ldots,\alpha_d] \in \Nset^d$, and denote the associated solution
of the PDE by $\uu_{\bm{\alpha}}$.

\section{Multi-Index Stochastic Collocation (MISC)}\label{section:MISC}

Besides IGA, the other basic building block of MISC is a tensorized quadrature
formula over the stochastic domain $\Gamma$ that can be used to
evaluate, e.g., expected value, variance and higher-order moments of
multi-variate random functions defined over $\Gamma$.

Let us therefore start by introducing a quadrature operator for
a univariate real-valued continuous function 
$v(t): \widetilde{\Gamma} \to \mathbb{R}$, with the understanding that $\widetilde{\Gamma}$
is a placeholder for any of the univariate sub-domains $\Gamma_1,\ldots,\Gamma_N$ composing the
random space $\Gamma = \Gamma_1 \times \ldots \times \Gamma_N$;
therefore, it comes with an associated probability density function,
that we denote by $\widetilde{\rho}$. The quadrature operator is defined as
\[
\mathcal{Q}^{m(\beta)}: C^0(\widetilde{\Gamma}) \to \mathbb{R}, \quad \mathcal{Q}^{m(\beta)}[v]=\sum^{m(\beta)}_{j=1} v(t_{\beta,j})\omega_{\beta,j},
\]
where $\beta \geq 1$ is a positive integer (usually referred to as the ``refinement level'' or just
``level'' of the quadrature operator),
$m(\beta)$ a strictly increasing function giving the number of distinct quadrature points to be used,
$\{ t_{\beta,j} \}_{j=1}^{m(\beta)} \subset \widetilde{\Gamma}$,
with corresponding weights $\{ \omega_{\beta,j} \}^{m(\beta)}_{j=1}$.
The quadrature points should be chosen according to the underlying probability measure $\widetilde{\rho}$;
see, e.g., \cite{nobile.eal:optimal-sparse-grids,ernst.eal:collocation-logn}.
Moreover, for refinement purposes it is advantageous if the quadrature points are chosen to be ``nested'', i.e., 
$\{ t_{\beta,j} \}_{j=1}^{m(\beta)} \subset \{ t_{\beta+1,j} \}_{j=1}^{m(\beta+1)}, \forall \beta \geq 1$.
In the numerical examples of this work we consider problems depending on uniform random variables,
for which a number of different families of nested quadrature points exist:
Clenshaw--Curtis points, several variants of Leja points, and Gauss--Patterson points;
see, e.g., \cite{nobile.eal:optimal-sparse-grids,trefethen:comparison,nobile.etal:leja,narayan:Leja,Chkifa:leja}.
In particular, we adopt the Clenshaw-Curtis points, which are defined as
\[
  t_{\beta,j}=\cos\left(\frac{(j-1) \pi}{m(\beta)-1}\right), \quad 1\leq j \leq m(\beta),
\]
and are nested provided that the function $m(\beta)$ is chosen as
\[
m(0)=0,\,\, m(1)=1,\,\, m(\beta)=2^{\beta-1}+1,\,\, \beta \geq 2.
\]


The extension to multi-variate real-valued continuous functions $v(\La{\yy}): \Gamma \to \mathbb{R}$
is obtained by tensorization of the univariate quadrature operators. In detail, we introduce a multi-index
$\bm{\beta} \in \Nset^N$, whose $i$-th component gives the level of the univariate quadrature to be
used along $\Gamma_i$, and we define the multi-variate quadrature operator as 
\begin{equation*}
  \mathcal{Q}^{m(\bm{\beta})}: C^{0}(\Gamma) \to \mathbb{R}, \quad
  \mathcal{Q}^{m(\bm{\beta})} = \bigotimes_{1 \leq i \leq N} \mathcal{Q}^{m(\beta_i)}, \quad
  \mathcal{Q}^{m(\bm{\beta})}[v]=\sum^{\# \bm{m}(\bm{\beta})}_{j=1} v(\La{\yy}_j)\omega_j,
\end{equation*}
where $\La{\yy}_j$ are points in the Cartesian grid $\bigotimes_{1 \leq i \leq N} \{z_{\beta_i,j}\}^{m(\beta_i)}_{j=1}$, $\omega_j$ are the corresponding products of weights
of the one-dimensional quadrature rules, and $\# \bm{m}(\bm{\beta})$ denotes the total number of
quadrature points in the Cartesian grid,  $\# \bm{m}(\bm{\beta}) = \prod_{i=1}^N m(\beta_n)$. 

We are now in the position to introduce the MISC approximation of the
expected value of a functional of the solution of the PDE, $\mathbb{E}[\Phi(\uu(\xx,\yy))]$.
For ease of presentation, we introduce the function $\phi:\Gamma \rightarrow \Rset$ which
associates each $\yy$ with its corresponding value of the functional $\phi(\yy)=\Phi(\uu(\xx,\yy))$;
thus, our goal becomes computing an approximation of $\mathbb{E}[\phi(\yy)]$.

Clearly, the value of $\phi(\yy)$ is accessible to us only by solving the PDE after fixing the
value of the random vector $\yy$: in particular, we denote by $\phi_{\bm{\alpha}}$ the value  
of $\phi$ obtained by post-process of $\uu_{\bm{\alpha}}$, where (as already mentioned) $\bm{\alpha}$
prescribes the number of elements in the knot vectors used to build the IGA solver, $N_{el,i} \sim 2^{\alpha_i}$.

The fully discrete approximation of $\mathbb{E}[\phi(\yy))]$ is therefore completely determined
by the choice of the discretizations $\bm{\alpha}$ in the physical space and by $\bm{\beta}$ in the probability space,
i.e.,
\[
  \La{\mathbb{E}[\phi(\yy))] \approx M_{\bm{\alpha},\bm{\beta}} = \mathcal{Q}^{\bm{m}(\bm{\beta})}[\phi_{\bm{\alpha}}].}
\]
Of course, the ideal approximation $M_{\bm{\alpha},\bm{\beta}}$ could be obtained by setting
$\alpha_1 = \ldots = \alpha_d = \bar{\alpha} \gg 1$ and $\beta_1=\ldots=\beta_N=\bar{\beta} \gg 1$, but this is out of reach
for even moderate values of $d,N,\bar{\alpha}$, and $\bar{\beta}$, due to its combinatorial computational cost.
Instead, we resort in MISC to the classical ``sparsification'' construction, which was already introduced in the literature
for solving high-dimensional PDEs \cite{b.griebel:acta} and quadrature problems \cite{smolyak:quadrature}.

To this end, we need to introduce the so-called univariate and multivariate ``detail operators''
on the physical and stochastic domains as follows,
with the understanding that $M_{\bm{\alpha},\bm{\beta}}=0$ when at least one component of $\bm{\alpha}$ or $\bm{\beta}$ is zero:
\begin{alignat*}{2}
&\text{\textbf{Univariate physical detail: }}
&& \Delta_i^{\text{phys}}[M_{\bm{\alpha},\bm{\beta}}]=M_{\bm{\alpha},\bm{\beta}}-M_{\bm{\alpha}-\bm{e}_i,\bm{\beta}} \text{ with } 1 \leq i \leq d; \\
&\text{\textbf{Univariate stochastic detail: }}
&& \Delta_i^{\text{stoc}}[M_{\bm{\alpha},\bm{\beta}}]=M_{\bm{\alpha},\bm{\beta}}-M_{\bm{\alpha},\bm{\beta}-\bm{e}_i} \text{ with } 1 \leq i \leq N; \\
&\text{\textbf{Multivariate physical detail: }}
&& \bm{\Delta}^{\text{phys}}[M_{\bm{\alpha},\bm{\beta}}] = \bigotimes_{i=1}^d \Delta_i^{\text{phys}}[M_{\bm{\alpha},\bm{\beta}}]; \\
&\text{\textbf{Multivariate stochastic detail: }}
&& \bm{\Delta}^{\text{stoc}}[M_{\bm{\alpha},\bm{\beta}}] = \bigotimes_{j=1}^N \Delta_j^{\text{stoc}}[M_{\bm{\alpha},\bm{\beta}}]; \\
&\text{\textbf{Mixed multivariate detail: }}
&& \bm{\Delta}^{\text{mix}}[M_{\bm{\alpha},\bm{\beta}}] = \bm{\Delta}^{\text{stoc}}\left[\bm{\Delta}^{\text{phys}}[M_{\bm{\alpha},\bm{\beta}}] \right].  
\end{alignat*}
Observe that taking tensor products of univariate details amounts to composing their actions, e.g.,
\[
\bm{\Delta}^{\text{phys}}[M_{\bm{\alpha},\bm{\beta}}]
= \bigotimes_{i=1}^d \Delta_i^{\text{phys}}[M_{\bm{\alpha},\bm{\beta}}]
= \Delta_1^{\text{phys}}\left[ \, \Delta_2^{\text{phys}}\left[ \, \cdots \Delta_d^{\text{phys}}\left[ M_{\bm{\alpha},\bm{\beta}} \right] \, \right] \, \right], 
\]
and analogously for the stochastic multivariate detail operators, $\bm{\Delta}^{\text{stoc}}[M_{\bm{\alpha},\bm{\beta}}]$.
Crucially, this in turn implies that the multivariate operators can be evaluated by evaluating certain full-tensor approximations
$M_{\bm{\alpha},\bm{\beta}}$ and then taking linear combinations:
\begin{align*}
\bm{\Delta}^{\text{phys}}[M_{\bm{\alpha},\bm{\beta}}]
& = \Delta_1^{\text{phys}}\left[ \, \Delta_2^{\text{phys}}\left[ \, \cdots \Delta_D^{\text{phys}}\left[ M_{\bm{\alpha},\bm{\beta}} \right] \, \right] \, \right] 
  = \sum_{\jj \in \{0,1\}^d} (-1)^{|\jj|} M_{\bm{\alpha}-\jj,\bm{\beta}};\\
\bm{\Delta}^{\text{stoc}}[M_{\bm{\alpha},\bm{\beta}}]
&  = \sum_{\jj \in \{0,1\}^N} (-1)^{|\jj|} M_{\bm{\alpha},\bm{\beta}-\jj}.
\end{align*}
The latter expression is known in the sparse-grids community as ``combination-technique'', and can be very
useful for practical implementations, especially for evaluating $\bm{\Delta}^{\text{phys}}[M_{\bm{\alpha},\bm{\beta}}]$:
indeed, it allows to evaluate detail operators by calling pre-existing
softwares on different meshes (in this case, IGA solvers) up to $2^d$ times in a ``black-box'' fashion.
As an example, if $d=2$, then
\begin{align*}
  \bm{\Delta^{\text{phys}}} [M_{\bm{\alpha},\bm{\beta}}]
  & = \Delta^{\text{phys}}_2 \left[ \, \Delta^{\text{phys}}_1 \left[ \, M_{\bm{\alpha},\bm{\beta}} \, \right] \, \right] \\
  & = \Delta^{\text{phys}}_2 [ M_{\bm{\alpha},\bm{\beta}} - M_{\bm{\alpha}-\ee_1,\bm{\beta}}] \\
  & = M_{\bm{\alpha},\bm{\beta}} - M_{\bm{\alpha}-\ee_1,\bm{\beta}} - M_{\bm{\alpha}-\ee_2,\bm{\beta}} + M_{\bm{\alpha}-\bm{1},\bm{\beta}}. 
\end{align*}
We remark that the four meshes needed to evaluate the combination-technique expression above are possibly anisotropic,
i.e., they may have different levels of discretization along the different physical directions.
A similar expression holds for the stochastic details, as well as for the mixed details.
Specifically, evaluating  $\bm{\Delta}^{\text{stoc}}[M_{\bm{\alpha},\bm{\beta}}]$
requires evaluating up to $2^N$ operators $M_{\bm{\alpha},\bm{\beta}}$ over different quadrature grids, and
evaluating $\bm{\Delta}^{\text{mix}}[M_{\bm{\alpha},\bm{\beta}}]$ requires
requires evaluating up to $2^{d+N}$ operators $M_{\bm{\alpha},\bm{\beta}}$ over different quadrature grids and physical meshes.
For instance, if $d=N=2$, then
\[
  \bm{\Delta^{\text{stoc}}} [M_{\bm{\alpha},\bm{\beta}}]
   = M_{\bm{\alpha},\bm{\beta}} - M_{\bm{\alpha}-\ee_1,\bm{\beta}} - M_{\bm{\alpha}-\ee_2,\bm{\beta}} + M_{\bm{\alpha}-\bm{1},\bm{\beta}};  
\]
\begin{align*}
  \bm{\Delta}^{\text{mix}}[M_{\bm{\alpha},\bm{\beta}}]
  & = \bm{\Delta^{\text{stoc}}} [\bm{\Delta^{\text{phys}}} [M_{\bm{\alpha},\bm{\beta}}] ] \\
  & = \bm{\Delta^{\text{stoc}}} [ M_{\bm{\alpha},\bm{\beta}} - M_{\bm{\alpha}-\ee_1,\bm{\beta}} - M_{\bm{\alpha}-\ee_2,\bm{\beta}} + M_{\bm{\alpha}-\bm{1},\bm{\beta}} ] \\
  & = M_{\bm{\alpha},\bm{\beta}} - M_{\bm{\alpha}-\ee_1,\bm{\beta}} - M_{\bm{\alpha}-\ee_2,\bm{\beta}} + M_{\bm{\alpha}-\bm{1},\bm{\beta}} \\
  &  - (M_{\bm{\alpha},\bm{\beta} - \ee_1} - M_{\bm{\alpha}-\ee_1,\bm{\beta} - \ee_1} - M_{\bm{\alpha}-\ee_2,\bm{\beta} - \ee_1} + M_{\bm{\alpha}-\bm{1},\bm{\beta}- \ee_1}) \\
  & - (M_{\bm{\alpha},\bm{\beta} - \ee_2} - M_{\bm{\alpha}-\ee_1,\bm{\beta} - \ee_2} - M_{\bm{\alpha}-\ee_2,\bm{\beta} - \ee_2} + M_{\bm{\alpha}-\bm{1},\bm{\beta} - \ee_2}) \\
  &  + M_{\bm{\alpha},\bm{\beta} - \bm{1}} - M_{\bm{\alpha}-\ee_1,\bm{\beta}  - \bm{1}} - M_{\bm{\alpha}-\ee_2,\bm{\beta}  - \bm{1}} + M_{\bm{\alpha}-\bm{1},\bm{\beta} - \bm{1}}.
\end{align*}
Observe that by introducing these operators we have access to a hierarchical decomposition of $M_{\bm{\alpha},\bm{\beta}}$;
indeed, the following is a telescopic identity
\[
M_{\bm{\alpha},\bm{\beta}} = \sum_{[\bm{\ii},\bm{\jj}] \leq [\bm{\alpha},\bm{\beta}]} \bm{\Delta}^{\text{mix}}[M_{\bm{\ii},\bm{\jj}}],
\]
i.e., it can be easily verified by replacing each term $\bm{\Delta}^{\text{mix}}[M_{\bm{\ii},\bm{\jj}}]$ with the
the corresponding combination-technique formula that all terms except $M_{\bm{\alpha},\bm{\beta}}$ will cancel.
For instance, recalling that by definition $M_{i,j} = 0$ when either $i=0$ or $j=0$, if $d=N=1$ we have
\begin{align}\label{eq:TP_example}
  \sum_{[i,i] \leq [2,2]} \bm{\Delta}^{\text{mix}}[M_{i,j}] 
  & = \bm{\Delta}^{\text{mix}}[M_{1,1}]
    + \bm{\Delta}^{\text{mix}}[M_{1,2}]
    + \bm{\Delta}^{\text{mix}}[M_{2,1}]
    + \bm{\Delta}^{\text{mix}}[M_{2,2}] \\
  & = M_{1,1}
    + ( M_{1,2} - M_{1,1} )
    + ( M_{2,1} - M_{1,1} )
    + M_{2,2} - M_{2,1} - M_{1,2} - M_{1,1} \nonumber \\
  & = M_{2,2}. \nonumber
\end{align}
The crucial observation is that (under suitable regularity assumptions on $\uu$),
not all of the details in the above hierarchical decomposition contribute equally to the approximation,
i.e., they can be discarded and the resulting formula will retain good approximation properties
at a fraction of the computational cost. Thus, we introduce the MISC approximation of $\mathbb{E}[\phi]$ as
\[
  \mathcal{I}^{\text{MISC}}_\Lambda=\sum_{[\bm{\alpha},\bm{\beta}] \in \Lambda} \bm{\Delta}^{\text{mix}}[M_{\bm{\alpha},\bm{\beta}}],
\]
for a suitable multi-index set $\Lambda \subset \mathbb{N}^{d+N}$,
which should be chosen as downward closed, see \eqref{eq:downward_closed_set}.
Clearly, the MISC estimator has a combination-technique expression as well, which can be written in compact form as
\begin{equation}\label{eq:misc_CT}
  \mathcal{I}^{\text{MISC}}_\Lambda
  =\sum_{[\bm{\alpha},\bm{\beta}] \in \Lambda} \bm{\Delta}^{\text{mix}}[M_{\bm{\alpha},\bm{\beta}}]
  =\sum_{[\bm{\alpha},\bm{\beta}] \in \Lambda} \sum_{\substack{[\bm{i},\bm{j}] \in \{0,1\}^{d+N}\\ [\bm{\alpha}+\bm{i},\bm{\beta}+\bm{j}] \in \Lambda}}
  (-1)^{\lvert [\bm{i},\bm{j}] \rvert} M_{\bm{\alpha},\bm{\beta}},
\end{equation}
which means that again we can evaluate MISC by evaluating full-tensor operators $M_{\bm{\alpha},\bm{\beta}}$ independently 
and combining them linearly according to  \eqref{eq:misc_CT}.

\begin{figure}
  \centering
  Isotropic full-tensor operator $M_{\bm{\alpha},\bm{\beta}}$\\
  \includegraphics[width=0.24\linewidth]{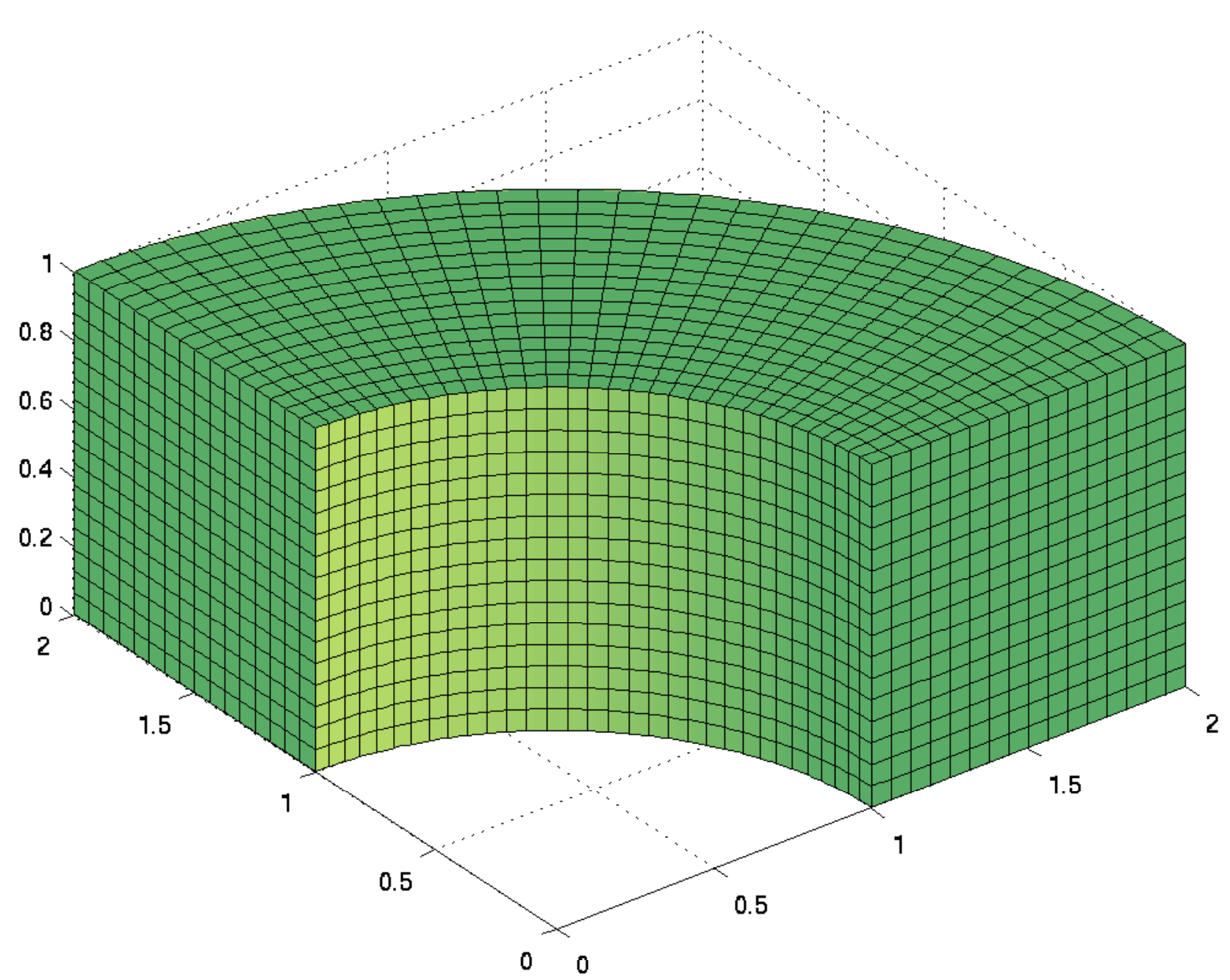} \quad \includegraphics[width=0.17\linewidth]{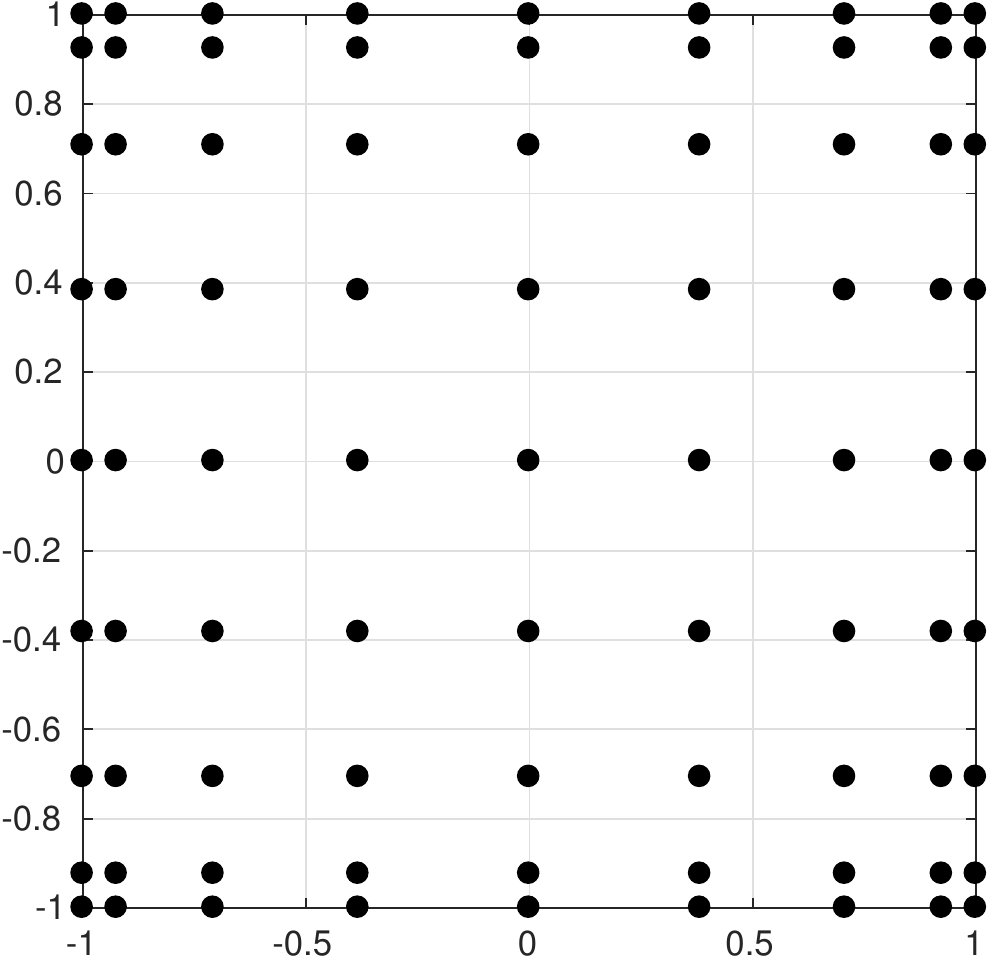}  \\[16pt]
  Anisotropic full-tensor operators $M_{\bm{\alpha},\bm{\beta}}$ for MISC \\[4pt]
  \includegraphics[width=0.24\linewidth]{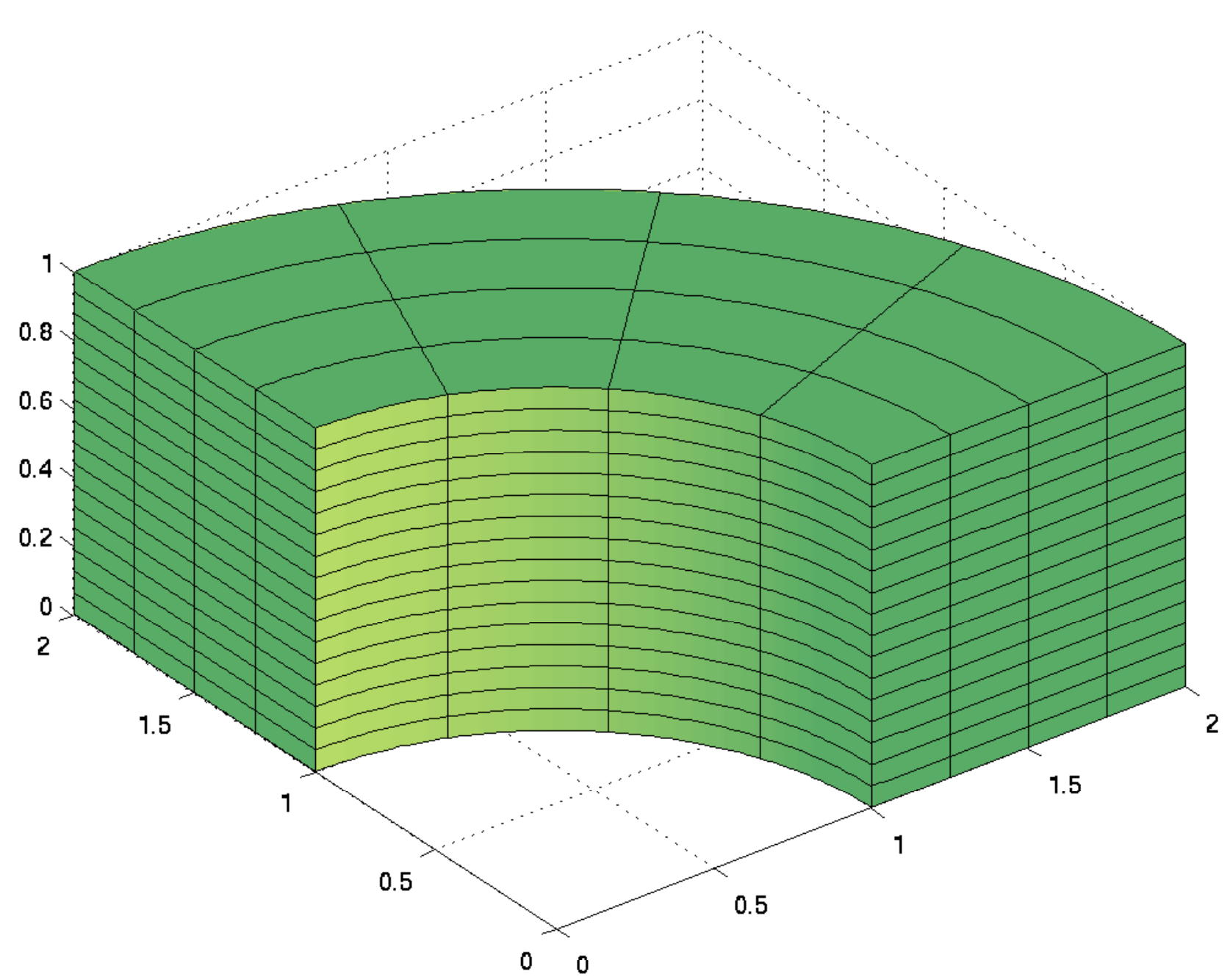} \quad \includegraphics[width=0.17\linewidth]{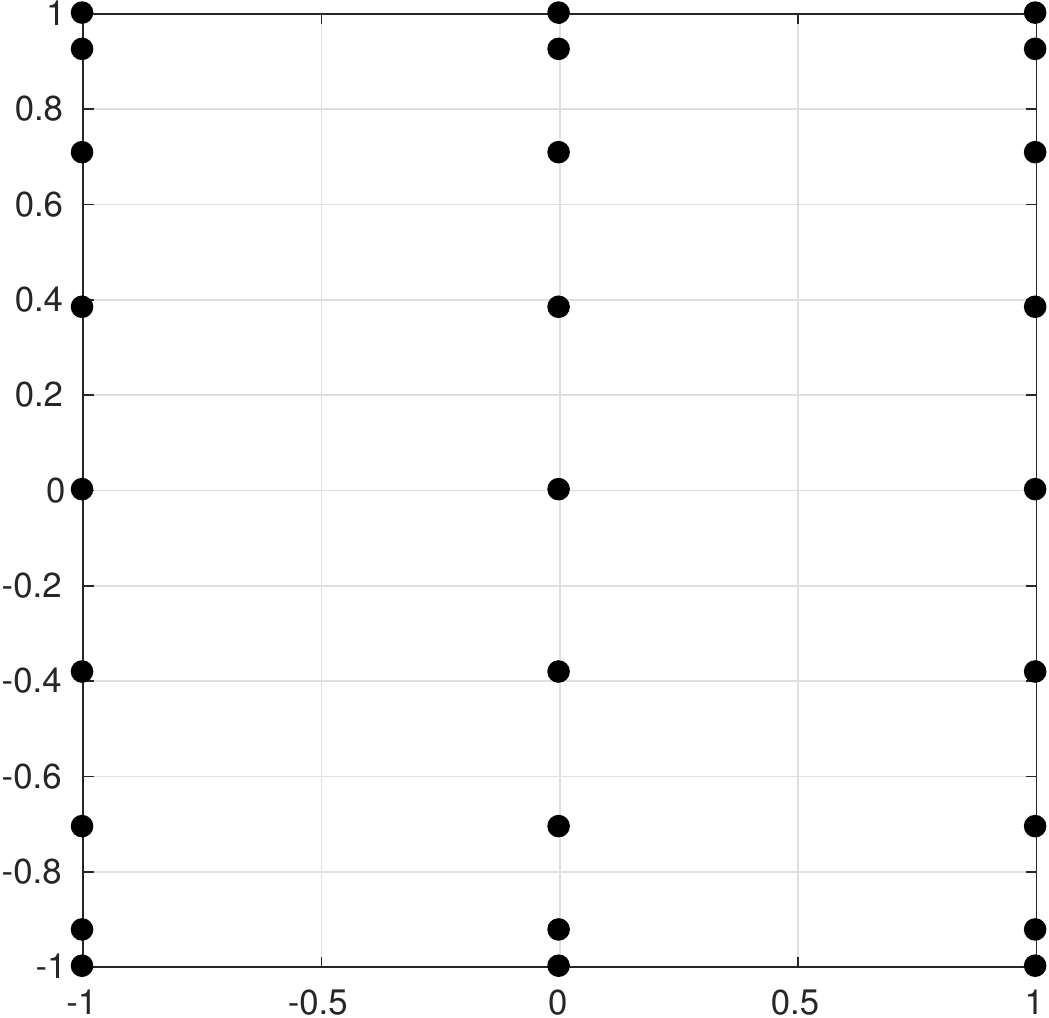}
  \qquad \qquad
  \includegraphics[width=0.24\linewidth]{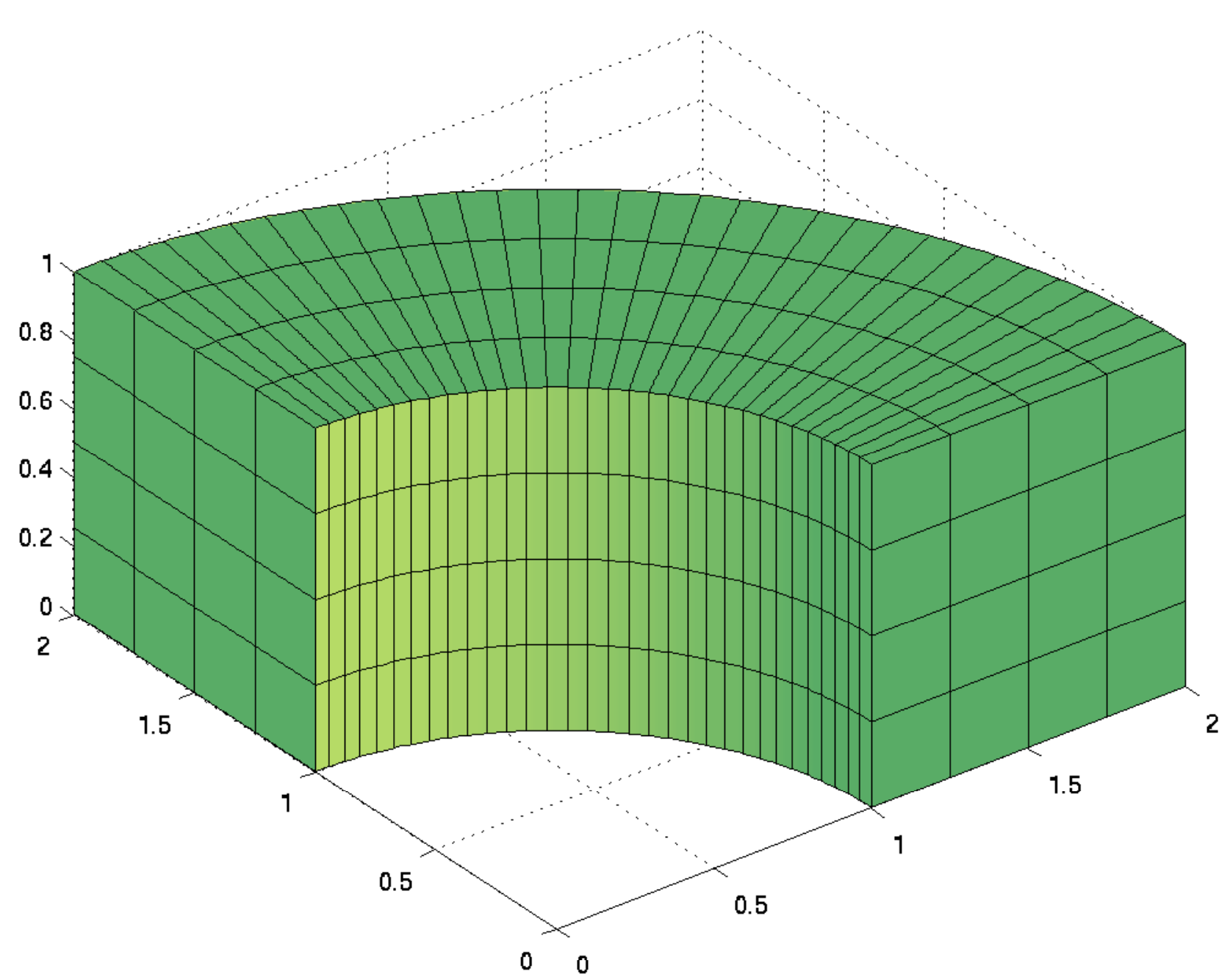} \quad \includegraphics[width=0.17\linewidth]{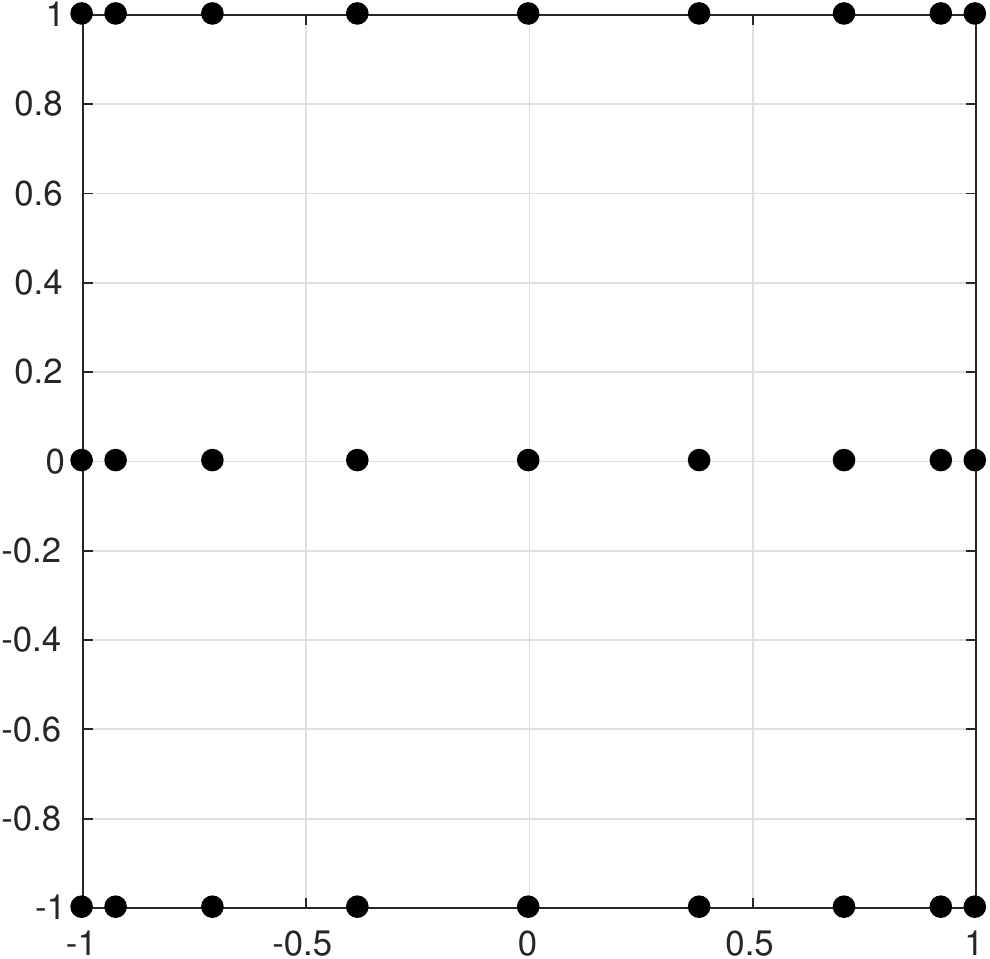}
  \caption{A full-tensor method requires sampling the stochastic domain with an isotropic Cartesian grid,
    and solving the PDE on a mesh refined along all physical directions for each sample, see top row.
    Instead, the MISC method requires that the stochastic domain be sampled with anisotropic Cartesian grids,
    and the PDE be solved on anisotropic meshes in space, which are combined in such a way that not all the physical
    and stochastic directions are simultaneously refined: some meshes satisfying this principle are shown in the
    second and third row. The results need to be linearly combined
    according the combination-technique formula in  \eqref{eq:misc_CT}.}\label{fig:MISC_in_picture}
\end{figure}

Of course, the effectiveness of the MISC estimator depends on the choice of the multi-index set $\Lambda$;
the general principle underlying the ``sparse'' construction is that $\Lambda$ should be chosen
to exclude isotropic full-tensor operators from the estimate, i.e., those operators that simultaneously
refine both the mesh on the physical domain and the quadrature grid on the stochastic domain.
Rather, indices that refine only a subset of the physical and/or stochastic directions should be chosen,
and then the combination-technique formula should be used to combine the partial results; this general idea is exemplified
in Figure \ref{fig:MISC_in_picture}. 
In a sense, the sparsification approach can be considered a Richardson extrapolation method, in the sense that
an ensemble of coarse approximations are combined in such a way that the result is more accurate than each of the components alone.
A simple yet quite effective choice of $\Lambda$ that keeps to a minimum the number of isotropic full-tensor operators to be computed
is the so-called ``total-degree'' set \cite{b.griebel:acta,back.nobile.eal:comparison}
\begin{equation}\label{eq:TD_MISC}
  \Lambda_{TD}(\bm{\kappa},\bm{g},w) = \left\{ [\bm{\alpha},\bm{\beta}] \in \Nset^{d+N} : \sum_{i=1}^d \kappa_i \alpha_i + \sum_{i=1}^N g_i \beta_i  \leq w  \right\},
  \text{ for some } w \in \Nset,
\end{equation}
where $\kappa_i,g_i$ are positive real values that can be used to allow more refinement along
selected physical or stochastic directions (the smaller the coefficient, the larger the maximum refinement level
allowed along that specific direction). 
For instance, if $d=N=1$, $w=2$, and the weights are chosen as $\kappa_i,g_i = 1$,
we have the set $\Lambda_{TD}(1,1,2) = \left\{ [i,j] \in \Nset^{2} : i + j  \leq 2  \right\}$, which amounts to
\begin{align*}
  \sum_{[i,j] \in \Lambda_{TD}(1,1,2)} \bm{\Delta}^{\text{mix}}[M_{i,j}] 
  & = \bm{\Delta}^{\text{mix}}[M_{1,1}]
    + \bm{\Delta}^{\text{mix}}[M_{1,2}]
    + \bm{\Delta}^{\text{mix}}[M_{2,1}] \nonumber \\
  & = M_{1,1}
    + ( M_{1,2} - M_{1,1} )
    + ( M_{2,1} - M_{1,1} ) \nonumber \\
  & = M_{1,2} + M_{2,1} - M_{1,1}. \nonumber  
\end{align*}
A similar choice of set $\Lambda$, i.e., total-degree sets with all weights set to 1,
was advocated also in \cite{bieri:sparse.tensor.coll,hps13}.
The last expression should be contrasted with \eqref{eq:TP_example}, and shows that the most refined
operator, $M_{2,2}$, no longer needs to be computed.
To devise an optimal strategy for selecting a good multi-index set, we introduce the error decomposition
\begin{align}
  \lvert \mathbb{E}[\phi] - \mathcal{I}^{\text{MISC}}_\Lambda \rvert
&  = \Big \lvert \mathbb{E}[\phi] - \sum_{[\bm{\alpha},\bm{\beta}] \in \Lambda} \bm{\Delta}^{\text{mix}}[M_{\bm{\alpha},\bm{\beta}}] \Big\rvert \nonumber \\
&  = \Big \lvert \sum_{[\bm{\alpha},\bm{\beta}] \not \in \Lambda} \bm{\Delta}^{\text{mix}}[M_{\bm{\alpha},\bm{\beta}}] \Big \rvert
  \leq \sum_{[\bm{\alpha},\bm{\beta}] \not \in \Lambda} \big \lvert \bm{\Delta}^{\text{mix}}[M_{\bm{\alpha},\bm{\beta}}] \big \rvert
  \leq \sum_{[\bm{\alpha},\bm{\beta}] \not \in \Lambda} E_{\bm{\alpha},\bm{\beta}},  \label{eq:err_decomp}
\end{align}
where we have defined $E_{\bm{\alpha},\bm{\beta}} = \big \lvert \bm{\Delta}^{\text{mix}}[M_{\bm{\alpha},\bm{\beta}}] \big \rvert$;
$E_{\bm{\alpha},\bm{\beta}}$ thus represents the ``error contribution'' of $[\bm{\alpha},\bm{\beta}]$,
i.e., the reduction in the approximation error due to having added $[\bm{\alpha},\bm{\beta}]$ to the current
index-set $\Lambda$. 
Similarly, we define the ``work contribution'' $W_{\bm{\alpha},\bm{\beta}}$ as the work required
to add $[\bm{\alpha},\bm{\beta}]$ to the current index-set $\Lambda$, for instance, summing the
degrees of freedom of all the new PDEs that need to be evaluated due to the addition of $[\bm{\alpha},\bm{\beta}]$. 
Thus, it can be easily seen that the strategy that delivers the best choice of $\Lambda$ consists of
adding to $\Lambda$ only the set of multi-indices with the largest profit
$P_{\bm{\alpha},\bm{\beta}} = \frac{E_{\bm{\alpha},\bm{\beta}}}{W_{\bm{\alpha},\bm{\beta}}}$
\cite{hajiali.eal:MISC1,hajiali.eal:MISC2}:
\begin{equation}\label{equation:set_of_largest_profits}
  \Lambda_\epsilon = \left\{ [\bm{\alpha},\bm{\beta}] \in \Nset^{d+N} : \frac{E_{\bm{\alpha},\bm{\beta}}}{W_{\bm{\alpha},\bm{\beta}}} \geq \epsilon \right\},
  \text{ for some } \epsilon > 0. 
\end{equation}
Such a set can be determined by classic adaptive algorithms such as those discussed
for quadrature problems in \cite{gerstner.griebel:adaptive,schillings.schwab:inverse,nobile.eal:adaptive-lognormal},
or according to a-priori bounds on the size of $E_{\bm{\alpha},\bm{\beta}}$, $W_{\bm{\alpha},\bm{\beta}}$;
in this work, we consider the latter approach. In particular, all the problems considered in the
numerical sections consist of elliptic PDEs, for which an expression of the optimal $\Lambda$ was derived in \cite{hajiali.eal:MISC1},
under the assumptions that the problem depends on uniform random variables and that the univariate quadrature operator
used on the stochastic domain is built over Clenshaw--Curtis points. In detail,
for some positive real values $r_i,c_i, i=1\ldots,d$ and $g_j, j=1\ldots,N$ that we will define in a moment,
the following bounds were introduced in \cite{hajiali.eal:MISC1}:
\begin{align}
  & E_{\bm{\alpha},\bm{\beta}} = \mathcal{O}(\widetilde{E}_{\bm{\alpha},\bm{\beta}}),
    \quad \widetilde{E}_{\bm{\alpha},\bm{\beta}} = 2^{-\sum_{i=1}^d \alpha_i r_i - \sum_{i=1}^N g_i 2^{\beta_i} \log_2 e}, \label{eq:E_estimate}\\
  & W_{\bm{\alpha},\bm{\beta}} = \mathcal{O}(\widetilde{W}_{\bm{\alpha},\bm{\beta}}),
    \quad \widetilde{W}_{\bm{\alpha},\bm{\beta}} = 2^{\sum_{i=1}^d \alpha_i c_i + \sum_{i=1}^N \beta_i}, \label{eq:W_estimate}
\end{align}
from which the following expression for (\ref{equation:set_of_largest_profits}) can be derived: 
\begin{equation}\label{eq:ABDO_MISC}
  \Lambda_{OP}(\bm{r},\bm{c},\bm{g},w) = \left\{ [\bm{\alpha},\bm{\beta}] \in \Nset^{d+N} :
      \sum_{i=1}^d (r_i+c_i) \alpha_i + \sum_{i=1}^N ( \beta_i +  g_i  2^{\beta_i} \log_2e) \leq w
  \right\}, \text{ for some } w \in \Nset,
\end{equation}
which is reminiscent of the total-degree set (\ref{eq:TD_MISC}).
The coefficients $r_i,c_i, i=1\ldots,d$ and $g_j, j=1\ldots,N$ are defined as follows,
and can either be determined a-priori or learnt during the execution of the MISC algorithm:
\begin{itemize}
\item $r_1,\ldots,r_d$ are the $h$-convergence rates of the PDE solver for the approximation of $\phi(\yy)$ for a fixed $\yy$, i.e.,
  \begin{equation}\label{eq:r_rate}
    |\phi_{\bm{\alpha}}(\yy) - \phi(\yy)| \leq C \prod_{i=1}^d N_{el,i}^{-r_i} = C \prod_{i=1}^d 2^{-\alpha_i r_i},
  \end{equation}
  where $\alpha_i$ can be determined a-priori by standard finite element theory, keeping into account
  possible corner/edge singularities as well as grading of meshes which may be introduced to
  mitigate the effects of such singularities; see, e.g., \cite{Guo1986,Babuska_book}.
  However, as already mentioned, we will measure $r_i$ numerically.
  \La{To this end, we perform some preliminary runs, in which we fix the stochastic parameters (to the midpoint of the parameter space)
    and increase the values of $\alpha_i$. We take the solution obtained with the largest value of $\alpha_i$ as a reference solution
    to replace  $\phi(\yy)$ in (\ref{eq:r_rate}), and compute $r_i$ e.g., by least square fitting. Note that we typically need only
    a few (4-5) of these preliminary runs for each $r_i$, and therefore the cost of this procedure is negligible
    with respect to the over-all cost of MISC. Moreover, we remark that while there might be a small dependence of $r_i$ on the specific
    value to which we fix $\yy$, this is typically not very severe, and we actually only need $r_i$ to be fixed to a reference value;
    in other words, the algorithm is typically robust to changes of the value of $\alpha$ from the second digit on.}

\item $c_1,\ldots,c_d$ are the rates of the increase of the cost of computing the approximation of $\phi(\yy)$ for a fixed $\yy$ as the physical mesh
  becomes finer,
  \begin{equation}\label{eq:c_rate}
    \text{cost}[\phi_{\bm{\alpha}}(\yy)] \leq C \prod_{i=1}^d 2^{\alpha_i c_i}.
  \end{equation}
  This cost is dominated by assembling and solving the linear system, since evaluating linear functionals
  of the solution is typically very cheap (e.g., a matrix-vector multiplication).
  We mention in-passing that relating the cost of assembling and solving
  the IGA linear system to the number of degrees of freedom
  is a delicate operation, see e.g. \cite{SERSH13}; we simply fit these rates from numerical experiments,
  \La{i.e., we just fit them with least squares, reusing the runs that were already performed to estimate $c_i$.}
  
\item $g_1,\ldots,g_N$ are the decay rates of the following bound on coefficients of the multivariate Legendre expansion of $\phi(\yy)$:
  \begin{equation}\label{eq:g_rate}
    \phi(\yy) = \sum_{\ii \in \Nset^N} \hat{\phi}_\ii \mathscr{L}_\ii(\yy),  \quad |\hat{\phi}_{\ii}|
    \leq C e^{-\sum_{n=1}^N g_n i_n}. 
  \end{equation}
  \La{We will also estimate numerically (by least square fitting) the coefficients $g_i$; moreover, note that
    we will fit directly the values of $\widetilde{E}_{\bm{\alpha},\bm{\beta}}$ in equation \eqref{eq:E_estimate}
    (see also Algorithm \ref{algo:misc_implementation}),
    so that in practice we do not need to compute the Legendre coefficients $\hat{\phi}_{\ii}$.}

\end{itemize}

In practice, we use the implementation reported in Algorithm \ref{algo:misc_implementation}.
The stopping criterion in Algorithm  \ref{algo:misc_implementation} above consists of checking
that the sum of the error contributions in the margin
of the multi-index set $\Lambda$, see (\ref{eq:margin}), is smaller than the required tolerance,
\begin{equation}\label{eq:MISC_stopping}
  \sum_{[\bm{\alpha},\bm{\beta}] \in \text{Mar}(\Lambda)} E_{\bm{\alpha},\bm{\beta}} < \text{TOL}.
\end{equation}
This choice of criterion is motivated by the error decomposition in (\ref{eq:err_decomp}), where we
further approximate the error bound by
\[
  |\mathbb{E}[\Phi(u)] - \mathcal{I}^{\text{MISC}}(u)|
  \leq \sum_{[\bm{\alpha},\bm{\beta}] \not \in \Lambda} E_{\bm{\alpha},\bm{\beta}}
  \approx \sum_{[\bm{\alpha},\bm{\beta}] \in \text{Mar}(\Lambda)} \widetilde{E}_{\bm{\alpha},\bm{\beta}}.
\]
The approximation above is only reasonable if the size of the details $E_{\bm{\alpha},\bm{\beta}}$
decreases quickly enough, which is true under sufficient smoothness hypotheses
for the problem at hand both with respect to the physical
variables and the random variables; see, e.g., \cite{Guignard:a-post}.

\begin{algorithm}[t]\label{algo:misc_implementation}
  \Fn{\AlCapSty{Multi Index Stochastic Collocation}($r_1,\ldots,r_d,c_1,\ldots,c_d,w_0,\text{TOL}$)}{
    $\Lambda = \Lambda_{TD}(\bm{1},\bm{1},w_0)$ \; 
    Compute MISC estimate $\mathcal{I}_{\Lambda}^{\text{MISC}}$ as in (\ref{eq:misc_CT}) \;
    Least squares fit of rates $g_i$ in (\ref{eq:E_estimate}) with $E_{\bm{\alpha},\bm{\beta}}$ for $[\bm{\alpha},\bm{\beta}] \in \Lambda$ \;
     Compute $\widetilde{E}_{\bm{\alpha},\bm{\beta}}$, $\widetilde{W}_{\bm{\alpha},\bm{\beta}}$ as in (\ref{eq:E_estimate}),  (\ref{eq:W_estimate}), and $\widetilde{P}_{\bm{\alpha},\bm{\beta}}=\widetilde{E}_{\bm{\alpha},\bm{\beta}}/\widetilde{W}_{\bm{\alpha},\bm{\beta}}$ for $[\bm{\alpha},\bm{\beta}] \in \text{Mar}(\Lambda)$ \;
    \While{$\sum_{[\bm{\alpha},\bm{\beta}] \in \text{Mar}(\Lambda)} \widetilde{E}_{\bm{\alpha},\bm{\beta}} > \text{TOL}$}{
    $\Theta=\left\{ [\bm{\alpha},\bm{\beta}] \in \text{Red}(\Lambda) : \widetilde{P}_{\bm{\alpha},\bm{\beta}}=\underset{[\bm{\alpha},\bm{\beta}] \in \text{Red}(\Lambda)}{\max}\{ \widetilde{P}_{\bm{\alpha},\bm{\beta}} \} \right\}$ \;
    $\Lambda=\Lambda \cup \Theta$ \;

     Compute MISC estimate $\mathcal{I}_{\Lambda}^{\text{MISC}}$ as in (\ref{eq:misc_CT}) \;
     Least squares fit of rates $g_i$ in (\ref{eq:E_estimate}) with $E_{\bm{\alpha},\bm{\beta}}$ for $[\bm{\alpha},\bm{\beta}] \in \Lambda$ \;
     Compute $\widetilde{E}_{\bm{\alpha},\bm{\beta}}$, $\widetilde{W}_{\bm{\alpha},\bm{\beta}}$ as in (\ref{eq:E_estimate}),  (\ref{eq:W_estimate}), and $\widetilde{P}_{\bm{\alpha},\bm{\beta}}=\widetilde{E}_{\bm{\alpha},\bm{\beta}}/\widetilde{W}_{\bm{\alpha},\bm{\beta}}$ for $[\bm{\alpha},\bm{\beta}] \in \text{Mar}(\Lambda)$ \;
    } 
  }
  \caption{MISC implementation}
\end{algorithm}

\begin{rmk}
  \La{In this work we have discussed an $h$-adaptive version of MISC, but in principle extending it to 
    a $p$ or even $hp$ version, in which the polynomial degree of the B-splines/NURBS IGA basis functions
    is also chosen adaptively, should not pose any conceptual challenge. One would need to extend the estimates
    of error and work contributions (\ref{eq:E_estimate}) and (\ref{eq:W_estimate})
    by adding an extra factor that takes into account the dependence of both
    quantities on $p$. The second ingredient would then be to introduce a proper ``quantization'' of $p$ refinement
    (i.e., whether one should increase the degree $p$ by 1 or more when refining)
    to make sure that the profit (i.e. error vs work ratio) obtained by raising the polynomial
    degree scales appropriately (i.e., that we are adding enough polynomials such that there is a significant reduction of the error when refining, and conversely that
    the work does not increase too much when refining). While to the best of our knowledge no paper proposes multi-level $hp$-adaptivity for UQ,
    a $p$-adaptive Multi-order Monte Carlo based on Discontinuous Galerkin solvers has been successfully proposed
    in \cite{mohamed:MOMC}.}
\end{rmk}

\section{Numerical results}\label{section:results}

In this section, we illustrate the performance of the MISC methodology using some numerical examples.
The random variables are considered to be uniformly distributed, and therefore we 
employ Clenshaw--Curtis quadrature points for the approximation over the stochastic space. 
Computational times were recorded on single-core runs of MISC on a workstation
equipped with Intel Xeon E5 processors with a clock rate of 2.8 GHz and an Ubuntu 16.04 operative system.
The IGA solver used is provided by the Matlab/Octave package GeoPDEs,
available at \url{http://rafavzqz.github.io/geopdes/}, see also \cite{VAZQUEZ2016523}.

\subsection{Test 1 - $3d$ linear elliptic PDE with random diffusion coefficient}

\begin{figure}[t]
  \centering
  \includegraphics[width=0.32\linewidth]{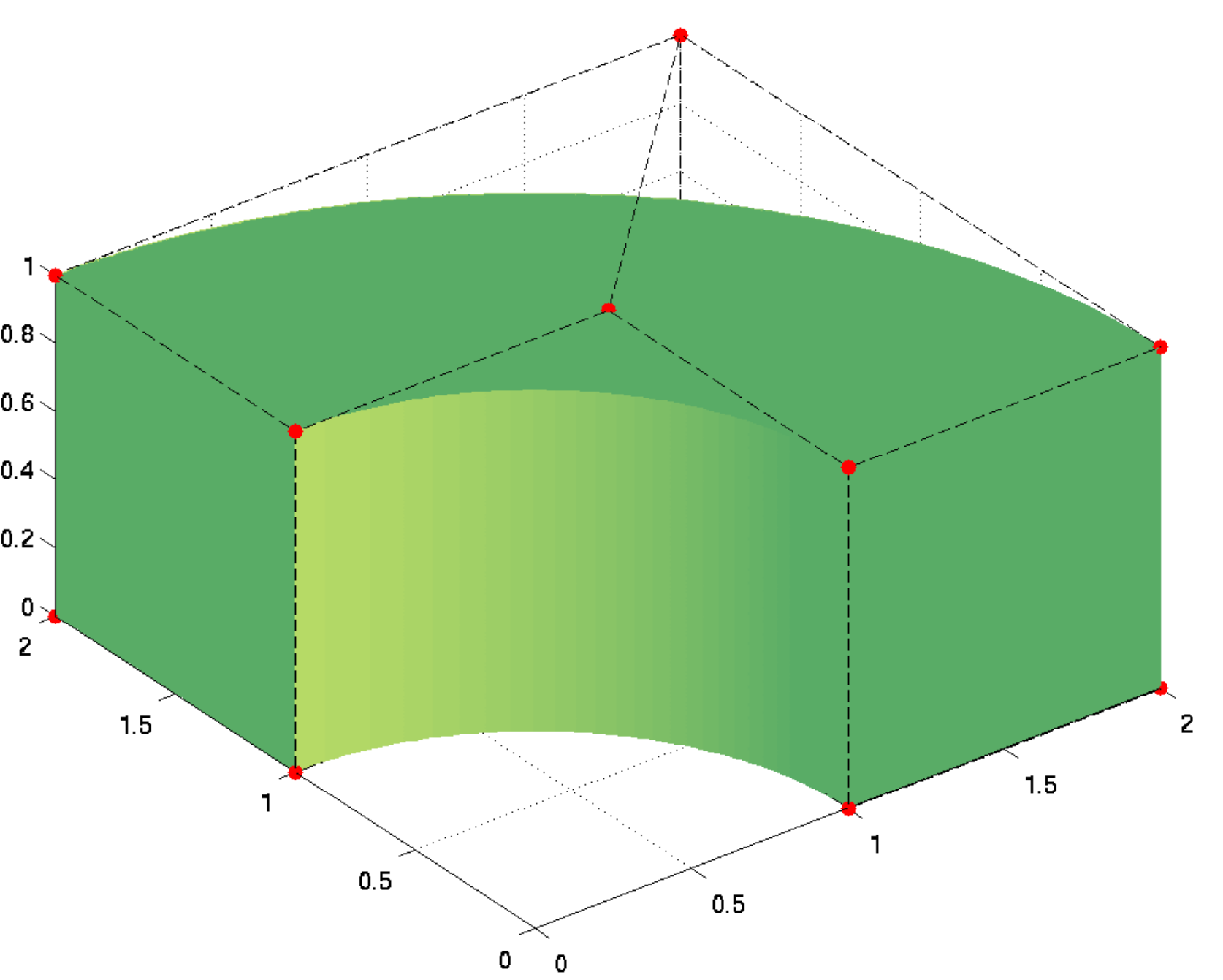}
  \includegraphics[width=0.32\linewidth]{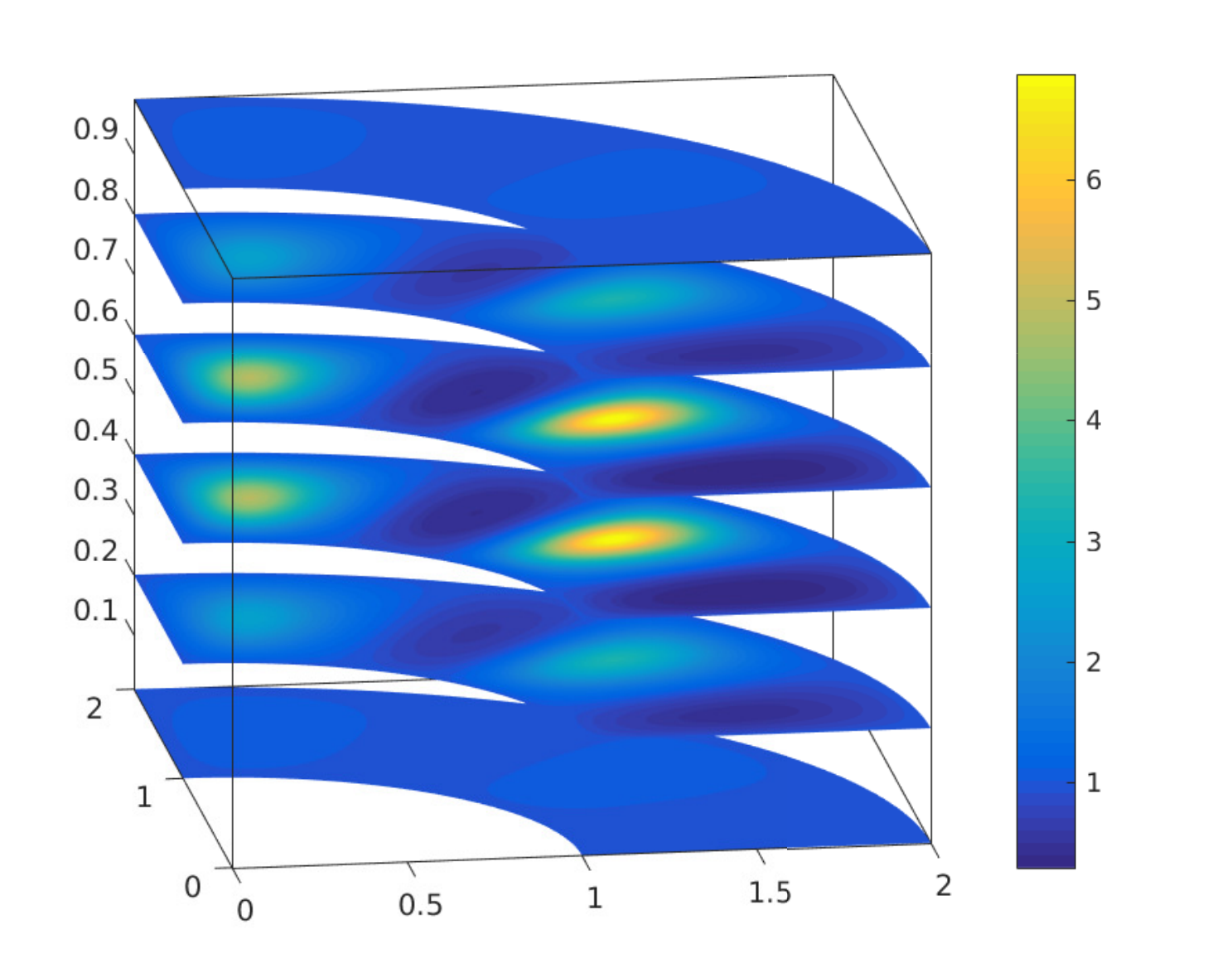}\\
  \includegraphics[width=0.32\linewidth]{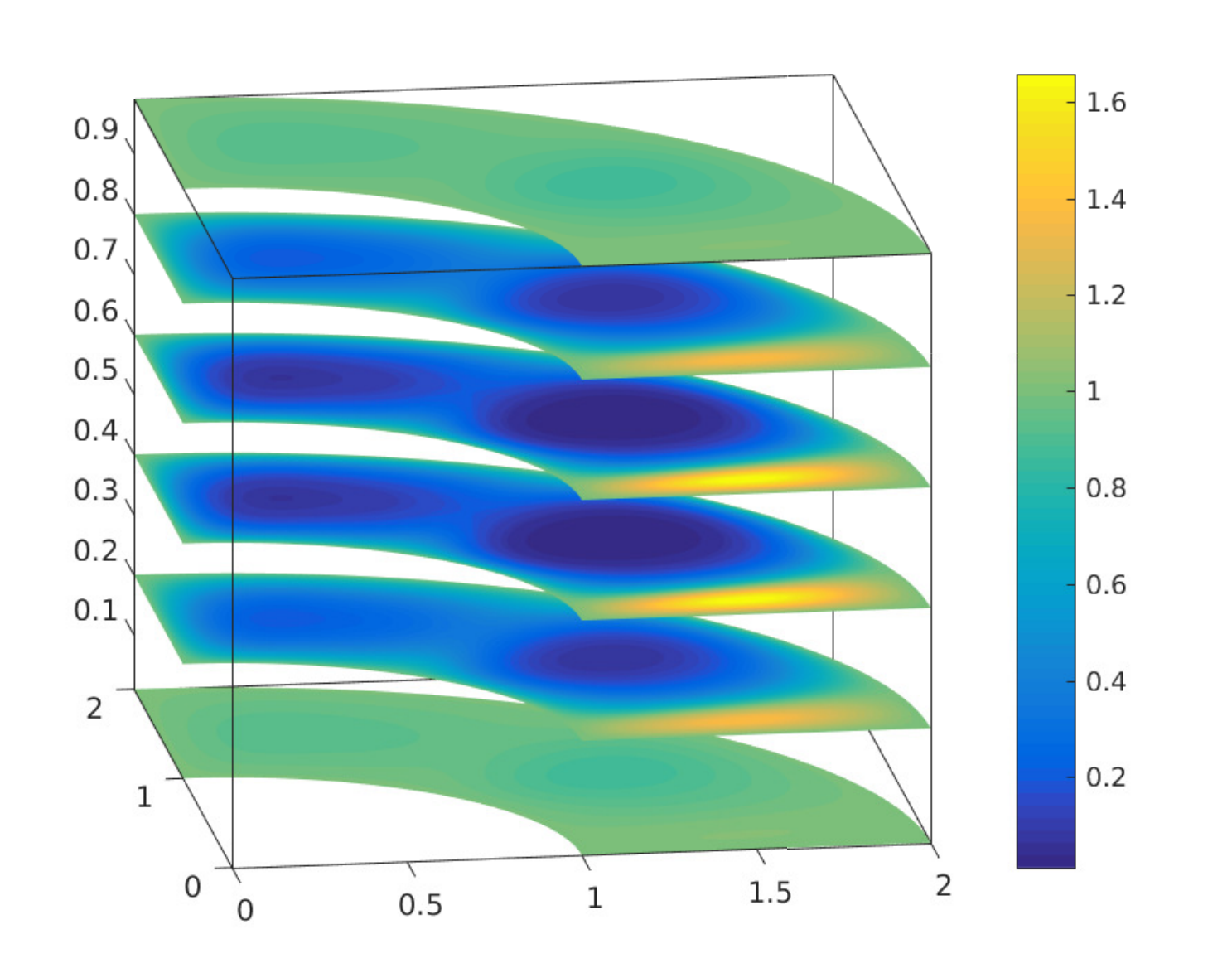}
  \includegraphics[width=0.32\linewidth]{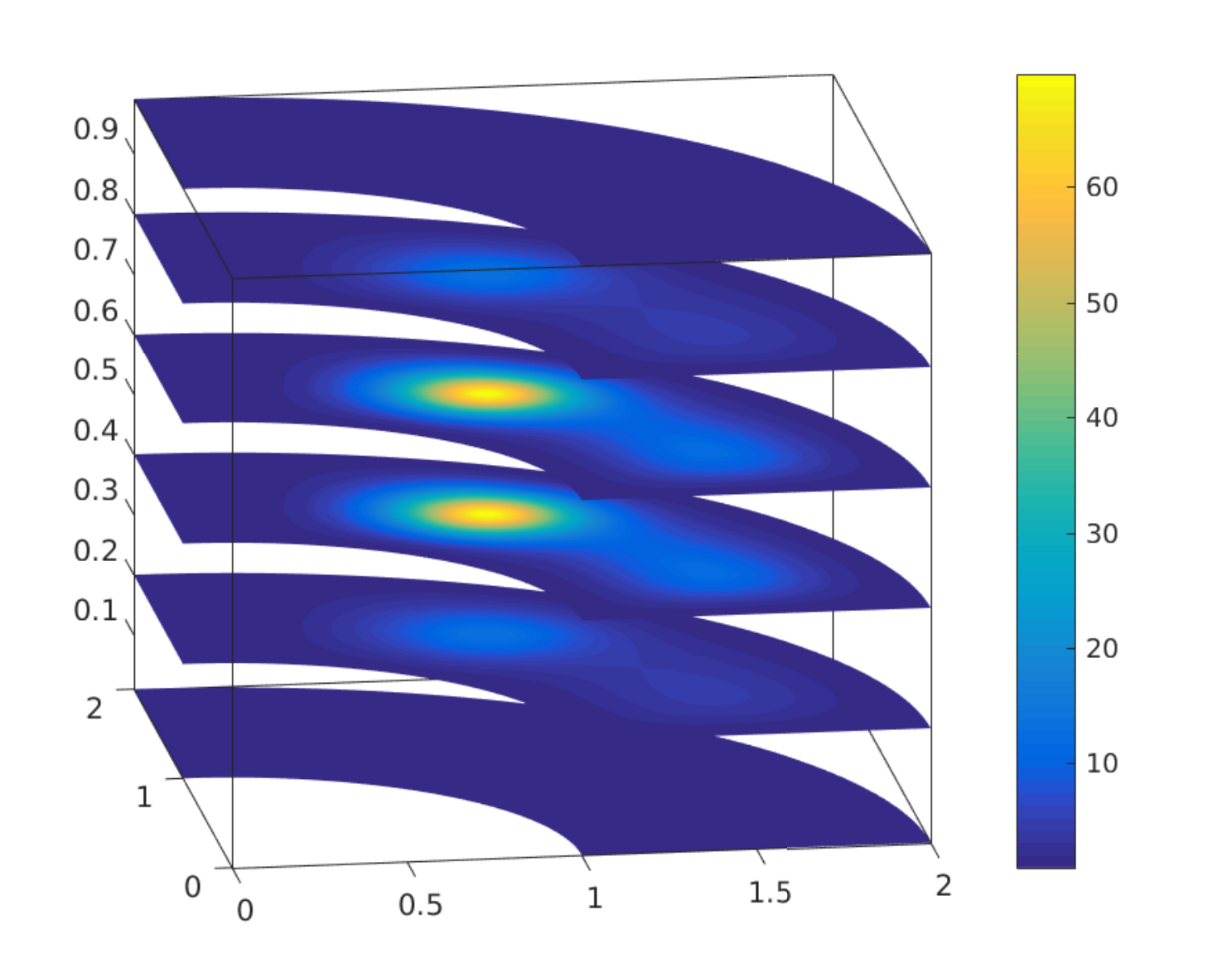}
  \caption{Test 1. Computational domain (top-left) and three different realizations of the random field.}\label{fig:test1_domain_and_randfield}
\end{figure}

In this test we consider a classic UQ benchmark, i.e., a linear elliptic PDE with random diffusion coefficient,
for which a vast body of literature exist, see e.g. \cite{babuska.nobile.eal:stochastic2,cohen_devore_2015} and references therein,
\[
  \begin{cases}
    -\text{div}[ a(\xx,\yy)\nabla u(\xx,\yy) ] = 1  & \xx\in \mcB,\\
    u(\xx,\yy)=0 & \xx\in\partial \mcB.
  \end{cases}
\]
We consider the ``thick quarter of ring'' in Figure \ref{fig:test1_domain_and_randfield}-top-left as physical domain;
this shape is a typical benchmark geometry in the IGA literature.
The random vector $\yy$ is composed of $N=3$ i.i.d. uniform random variables,
$y_i \sim \mathcal{U}(-1,1)$, i.e., $\Gamma = [-1,1]^3$. The random field $a(\xx,\yy)$ models the variability
in the properties of the material, e.g., uneven heat capacity due to imperfections.
Because of the peculiar shape of the computational domain, we express the random field in cylindrical coordinates: 
\begin{alignat*}{2}
  &  a([\rho, \theta,z],\yy) =      && e^{c\, \gamma([\rho, \theta,z],\yy)}, \quad \text{with } c=4 \text{ and}\\
  & \gamma([\rho, \theta,z],\yy) =  &&      y_1 \sin\left(2 \theta \right) \sin\left( \pi(\rho-1) \right) \sin\left(\pi z\right) + \\
  &                                 && 0.4  y_2 \sin\left(8 \theta \right) \sin\left( \pi(\rho-1) \right) \sin\left(\pi z\right) + \\
  &                                 && 0.1  y_3 \sin\left(16 \theta \right) \sin\left( \pi(\rho-1) \right) \sin\left(\pi z\right).
\end{alignat*}
The expression for $\gamma([\rho, \theta,z],\yy)$ mimics the expression that one would obtain
by applying a spectral decomposition like Fourier \cite{back.nobile.eal:lognormal} or Karhunen--Lo\`eve \cite{ghanem.spanos:book}
to a random field and then truncating it to retain only the most important modes.
Three different realizations of the random field can be seen in Figure \ref{fig:test1_domain_and_randfield}.
Note the different scales of the point-wise values of the realizations (due to the magnifying effect of the exponential
operation in the definition of the random field), as well as the differences among the frequencies of oscillations.
By construction, there exist two real values $0<m<M$ such that $ m < a([\rho, \theta,z],\yy) < M$
for $\rho$-almost every $\yy \in \Gamma$, which guarantees that the problem is well-posed in $V=H^1_0(\mcB)$,
and that there exists an optimally convergent approximation of $u$ in $L^2_\rho(\Gamma;V)$ as well as $L^\infty_\rho(\Gamma;V)$
based on either $\rho$-orthonormal polynomials or interpolation processes; see, e.g., \cite{nobile.eal:optimal-sparse-grids,chkifa:nonlinear}.
We are interested in computing the expected value of the integral of the solution over the physical domain,
$\Phi(\vv) = \int_\mcB \vv(\xx,\cdot)d\xx, \forall \vv \in V$.

Concerning the IGA solver, note that the computational domain in this example cannot be described exactly by
B-splines, because some of the edges are circle arcs. Therefore, as already pointed out in Section \ref{section:IGA},
we resort to using NURBS basis functions instead of B-splines, which however does not imply any change in the MISC
procedure. Specifically, we employ NURBS of degree $p=2$ and maximal continuity, $C^1$;
incidentally, note that this means that a degree-elevation operation (i.e., $p$-refinement) will have to be performed prior
to actually starting the MISC computation, since the ``thick quarter of ring'' geometry is defined
by linear polynomials in two out of three parametric directions.
The knots in the parametric domain are not uniform; rather, they are scaled towards the edges of the domain
according to a power law with exponent set to 3 to capture the edge singularities of the PDE solution for fixed $\yy$,
and thus improve the convergence of the solution as the mesh-size decreases.
We chose the exponent based on a numerical exploration
aimed at recovering the optimal convergence of the IGA solver on analogous problems
(albeit without randomness), see \cite{beck.eal:sparse-IGA} and references therein for more details.
For a fixed realization of the random field, the Galerkin stiffness matrix is typically
less sparse for B-splines and NURBS than for finite elements (due to the larger support of B-splines/NURBS,
which is proportional to the degree $p$, \cite{SERSH13}), so we employ a direct solver (Matlab's backslash).

We will compare the convergence results of MISC with two different multi-level methods, namely,
the Multi-Level Monte Carlo (MLMC, \cite{scheichl.giles:MLMC}) in the implementation proposed in
\cite{hajiali.eal:continuationMLMC}, and its refined version
Multi-Index Monte Carlo (MIMC, \cite{hajiali.eal:MultiIndexMC}). 

The convergence of MLMC and MIMC for elliptic PDEs with random coefficients
was discussed respectively in \cite{scheichl.giles:MLMC} and in \cite{hajiali.eal:MultiIndexMC}
and depends on the rates $r_i, c_i$ in \eqref{eq:r_rate} and \eqref{eq:c_rate}. 
We have found numerically that $r_i, c_i$ have approximate values of $r_i=[4,4,4]$ and $c_i=[1,1,1]$,
which implies that the computational cost for reaching an accuracy $\text{TOL}$ 
is expected to be $\mathcal{O}(\text{TOL}^{-2})$ for both MLMC and MIMC,
see \cite{scheichl.giles:MLMC,hajiali.eal:MultiIndexMC} for details.
This is the optimal rate for sampling schemes, i.e., most of the sampling is done
on the coarsest mesh levels, such that the computational cost is ``equivalent
to sampling a random variable'' (of course, up to the cost of solving the linear
system corresponding to the coarsest levels, which is still non-negligible). 
The cost of a standard Monte Carlo analysis where, for an assigned tolerance $\text{TOL}$,
we choose a physical mesh and a number of samples in the stochastic domain such
that both the deterministic and statistical errors are smaller than $\text{TOL}/2$,
would instead be proportional to $\mathcal{O}(\text{TOL}^{-2.75})$; see, e.g., \cite{scheichl.giles:MLMC}.

Concerning MISC, convergence results are available in \cite{hajiali.eal:MISC1} for
the case of a PDE depending on a finite number of random variables (as is the case in the current example)
and in \cite{hajiali.eal:MISC2} for the case of a PDE depending on a countable sequence of random variables.
The convergence result for finitely many random variables in \cite{hajiali.eal:MISC1}
depends on the rates $r_i, c_i$ only, i.e., $g_i$ do not play any role;
in other words, the approximation over the probability space by tensorized quadrature
is expected to converge fast enough that it would not impact the overall convergence rate.
More precisely, the result in \cite{hajiali.eal:MISC1} predicts for the a-priori chosen set in (\ref{eq:ABDO_MISC})
the asymptotic convergence estimate 
\[
  |\mathbb{E}[\Phi(u)] - \mathcal{I}^{\text{MISC}}(u)| \leq C \text{Work}^{-4} (\log \text{Work})^{10},
\]
which would in turn imply that the computational cost of reaching an accuracy $\text{TOL}$ 
is expected to be asymptotically $\mathcal{O}(\text{TOL}^{-1/4})$ up to logarithmic terms.

\begin{figure}[t]
  \centering
  \includegraphics[width=0.42\linewidth]{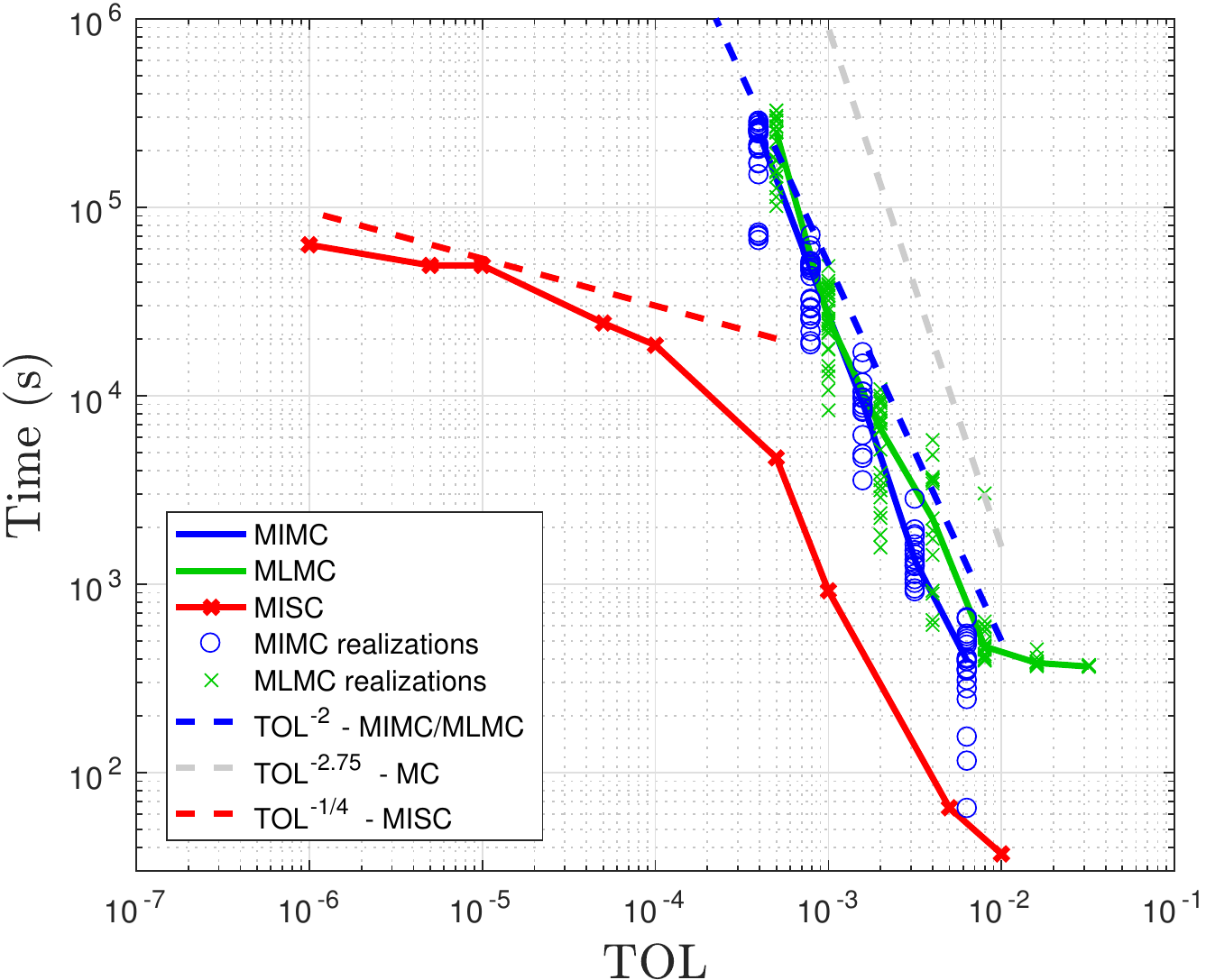}\qquad
  \includegraphics[width=0.42\linewidth]{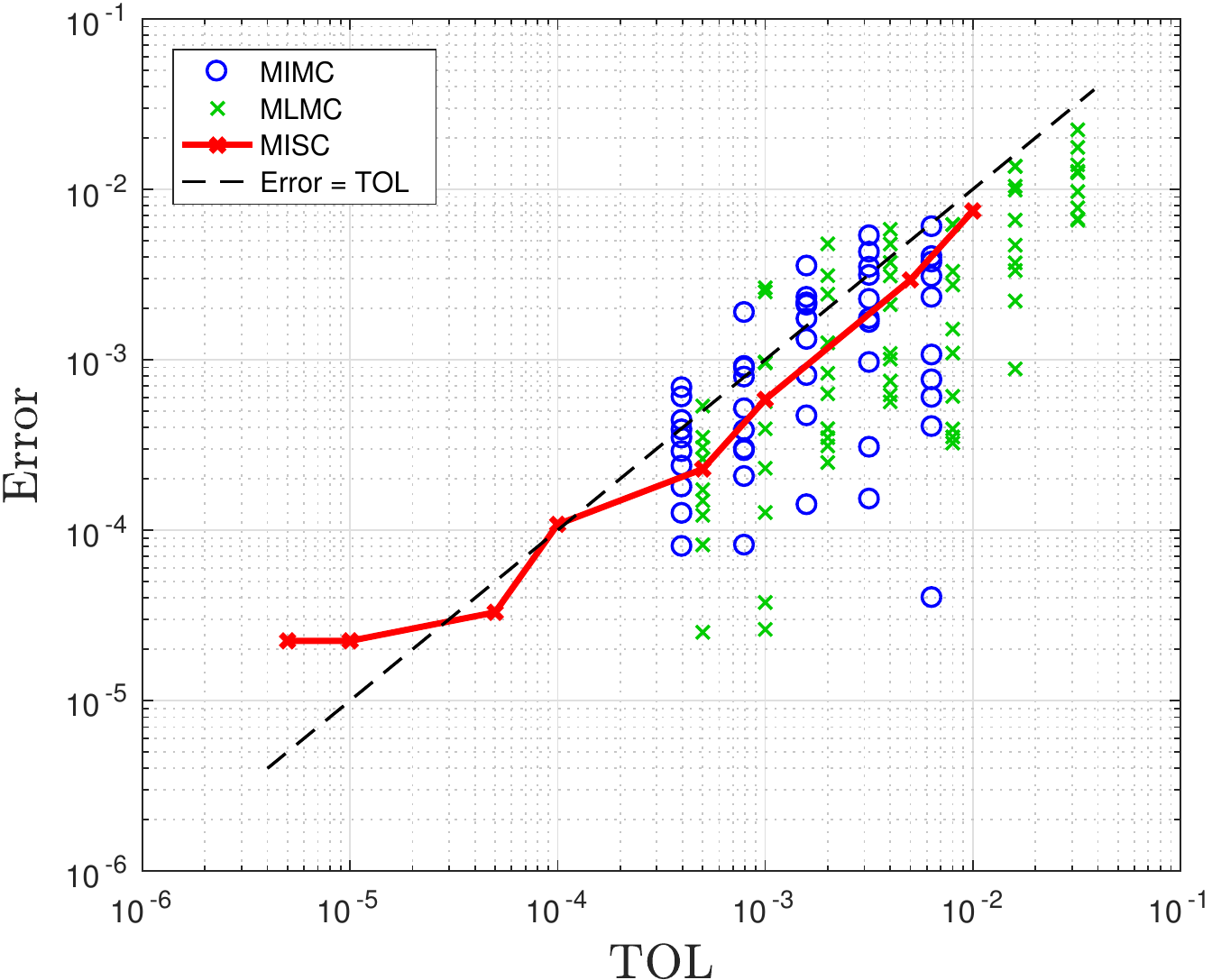}
  \caption{Test 1. Left: Convergence results. Right: Consistency of stopping criterion; the error has been computed against a
    sufficiently refined MISC solution. Markers for MLMC and MIMC show the value of the error attained for each method run.
    We have enforced that the error for MLMC and MIMC should be less than TOL with 95\% probability (asymptotically).}\label{fig:test1_results}
\end{figure}

We report the computational results in Figure \ref{fig:test1_results}-left. 
On the horizontal axis we show the tolerance, TOL, which is used in the stopping criterion of each algorithm,
and on the vertical axis the recorded computational time:
thus, the flatter the convergence curve, the more effective the method, i.e., moving to smaller tolerance
does not require a dramatic increase in the computational time.
Of course, the implicit assumption behind this plot is that the error actually achieved
once the algorithm stops is similar in size to the tolerance enforced as 
stopping criterion. For MLMC and MIMC, this is guaranteed with ``high probability''
by the choice of number of samples per mesh, while for MISC we employ the criterion in (\ref{eq:MISC_stopping});
see Figure \ref{fig:test1_results}-right for the effectiveness of these stopping criteria.

We observe from the numerical results that the convergence rates of MLMC and MIMC are roughly
2, in agreement with the theory discussed above. Observe also that the performances of MLMC and MIMC are very
close, due to the fact that the convergence of the IGA approximation of the physical problem
is very fast ($r_i=4$), hence sampling the physical problem over anisotropic meshes in space does not provide as significant advantage.
MISC features instead a convergence rate of roughly $1/4$, as discussed above, after a substantial
pre-asymptotic regime which was also expected, \La{due to the fact that the diffusion field features
  quick oscillations which cannot be properly captured by the coarse grids employed at the beginning of
  the simulation}. Moreover, MISC has a significantly smaller error, which testifies its superior performance compared
to the other methods reported, in agreement with previous numerical investigations
reported in \cite{hajiali.eal:MISC1,hajiali.eal:MISC2}.

\begin{figure}[t]
   \centering
   \includegraphics[width=0.42\linewidth]{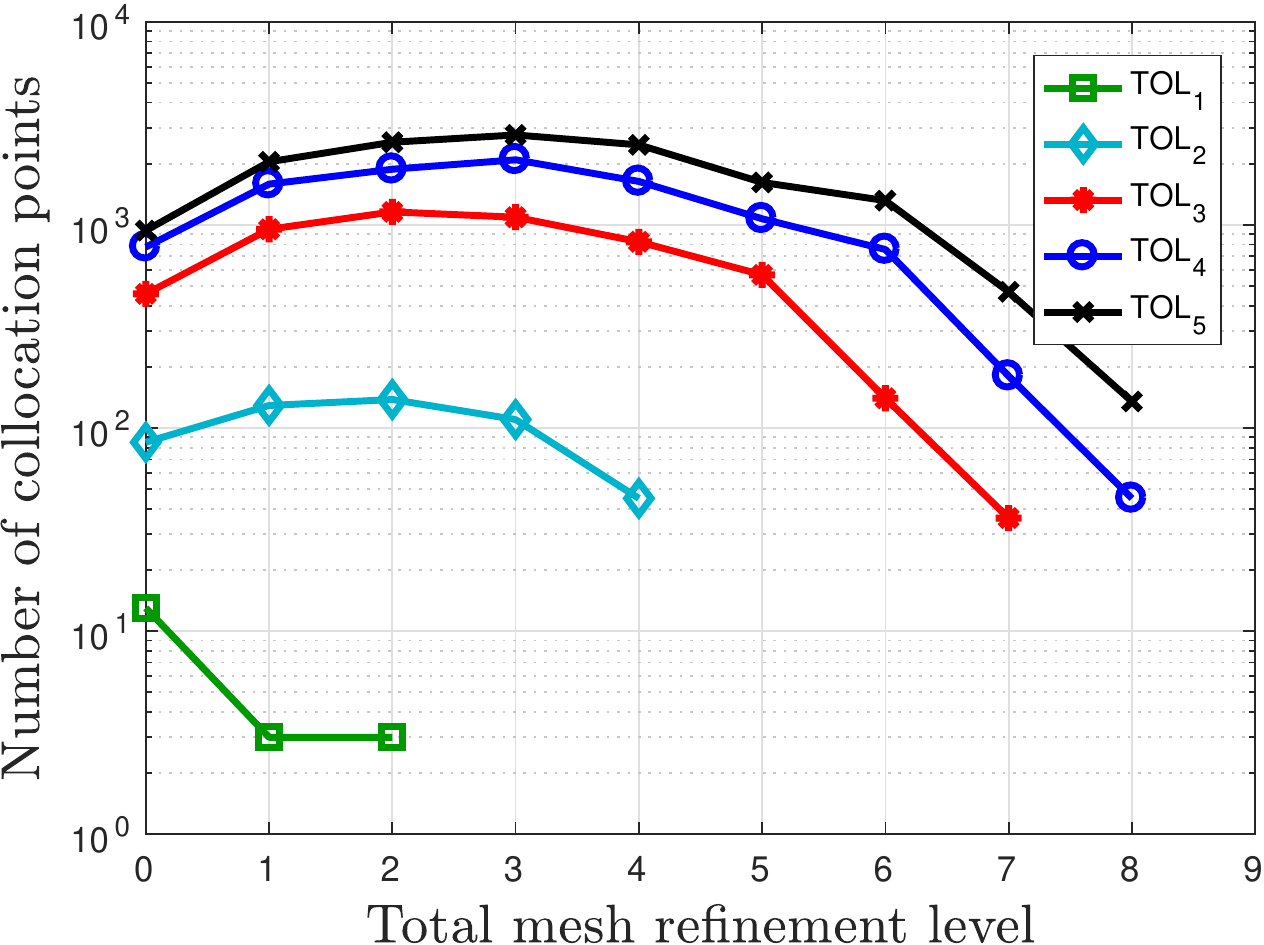}
   \caption{\La{MISC computational cost allocation.}}
   \label{fig:cost}
 \end{figure}

\begin{figure}[t]
   \centering
   \includegraphics[width=0.46\linewidth]{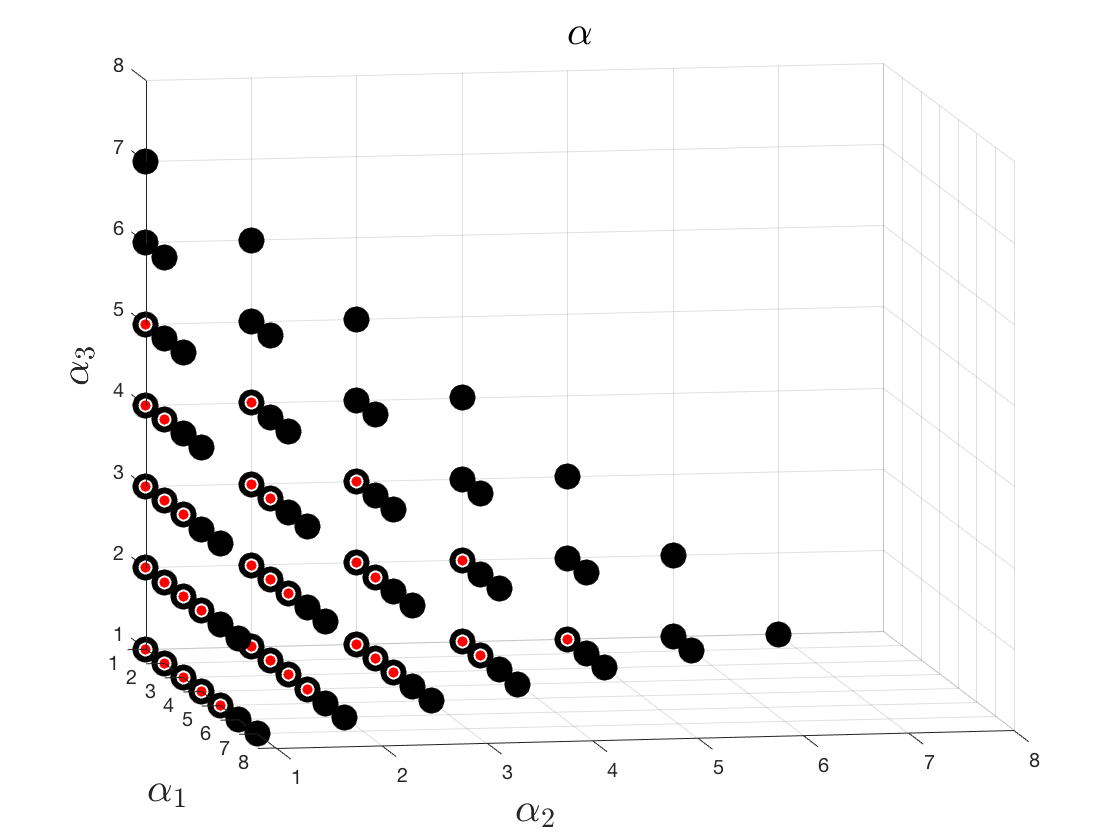}\qquad \quad
   \includegraphics[width=0.46\linewidth]{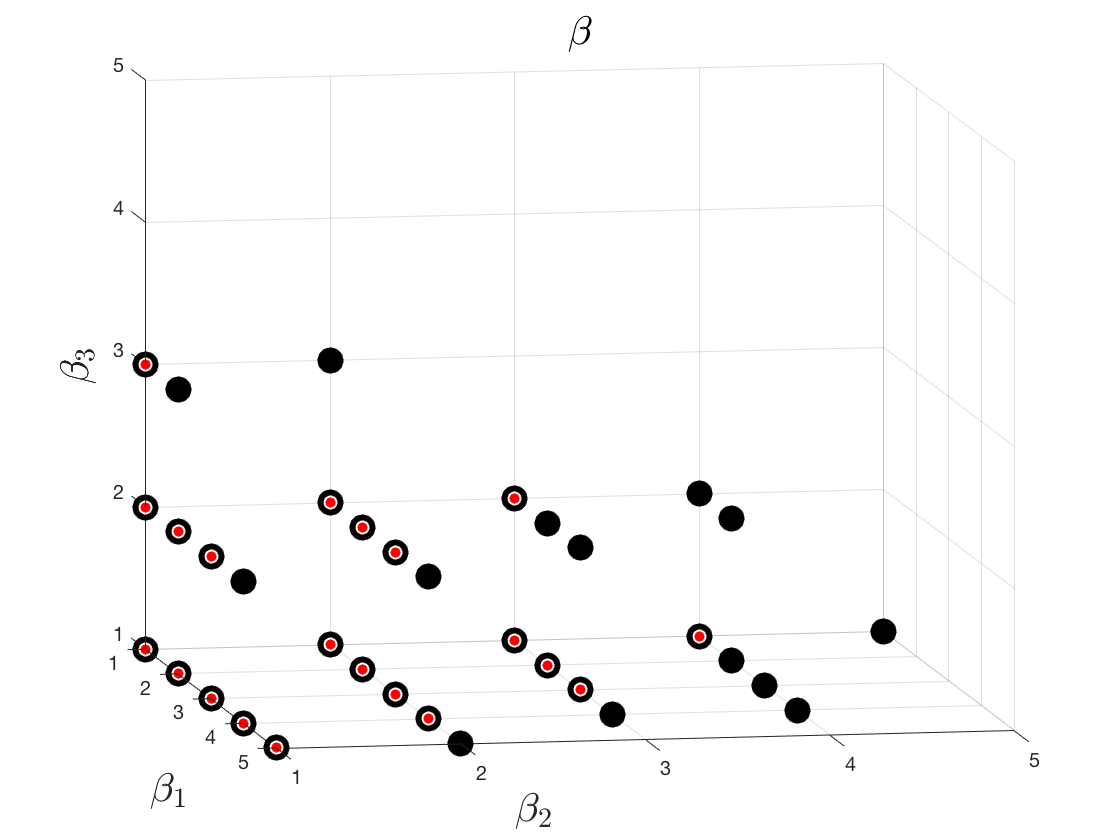}\\[12pt]
   \includegraphics[width=1.0\linewidth]{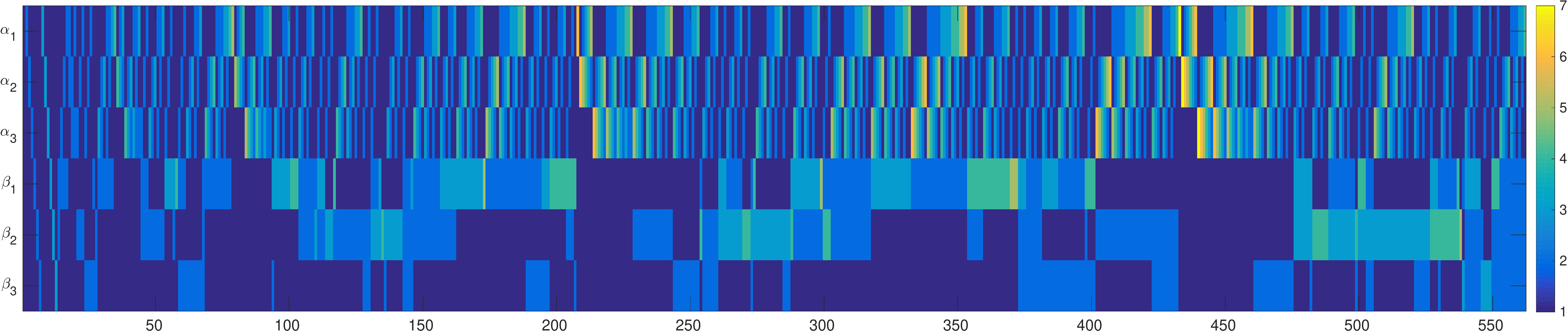}
   \caption{\La{The top figures show $\bm{\alpha}$-indices (top-left) and $\bm{\beta}$-indices (top-right) of the MISC multi-index set $\Lambda$ for $\hbox{TOL}=10^{-3}$ (red dots) and $\hbox{TOL}=10^{-4}$ (black dots), respectively. The bottom figure shows the multi-index values, $[\bm{\alpha},\bm{\beta}]$, in the order they were added to $\Lambda$. Each multi-index is shown as a segmented vertical line, colored ranging from blue to yellow according to the values of each entry of the multi-index.}}
   \label{fig:indexset}
 \end{figure}

\La{In order to study how the MISC method allocates the computational effort, let us now define \emph{total mesh refinement} as $\sum^{d}_{i=1} (\alpha_i-1)$. Note that, under the given work model \eqref{eq:c_rate} and with $c_i$ being identical in all directions (as in our case), the physical meshes of the same total mesh refinement level have also the same computational work. In Figure \ref{fig:cost}, we observe that the number of collocation points allocated to the expensive meshes decreases fast with the total mesh refinement level. Interestingly enough, the coarsest physical meshes are not guaranteed to be the most evaluated, which is in contrast to the standard MLMC. 
  The top plots of Figure \ref{fig:indexset} show the multi-indices $\bm{\alpha}$ and $\bm{\beta}$ in $\Lambda$ for $\hbox{TOL}=10^{-3}$ and $\hbox{TOL}=10^{-4}$. We see that the set of $\bm{\alpha}$-indices in $\Lambda$ is a total-degree set, see \eqref{eq:TD_MISC}, with maximum refinement level $5$ and $7$ for $\hbox{TOL}=10^{-3}$ and $\hbox{TOL}=10^{-4}$, respectively, whereas the set of $\bm{\beta}$-indices shows that the $\beta_3$ direction requires less refinement. The bottom plot in Figure \ref{fig:indexset} shows the multi-index values, $[\bm{\alpha},\bm{\beta}]$, in the order they were added to $\Lambda$, and here we observe the sparsity of the multi-index set $\Lambda$
  as no single multi-index has high values in more than a few components, in contrast to its full-tensor counterpart; in other words, the vertical
  lines (each showing the value of the components of a multi-index) are predominantly showing low values (shades of blue), instead of large values 
  (green to yellow shades).}

\La{Finally, we analyze the convergence of the estimates of the rates $g_i$, see equation \eqref{eq:g_rate}. These are estimated progressively by the MISC method, as shown in  Algorithm \ref{algo:misc_implementation}, and the estimated values are indeed converging as we keep adding more multi-indices to $\Lambda$, or in other words, as we decrease the required tolerance $\hbox{TOL}$. For Test 1, as expected, the parameter $g_3$ has the largest value, which is consistent with the lower maximum refinement observed in direction $\beta_3$ of the $\bm{\beta}$-indices set, as shown in Figure \ref{fig:indexset}-top-right.}

 \begin{figure}[t]
   \centering
   \includegraphics[width=0.42\linewidth]{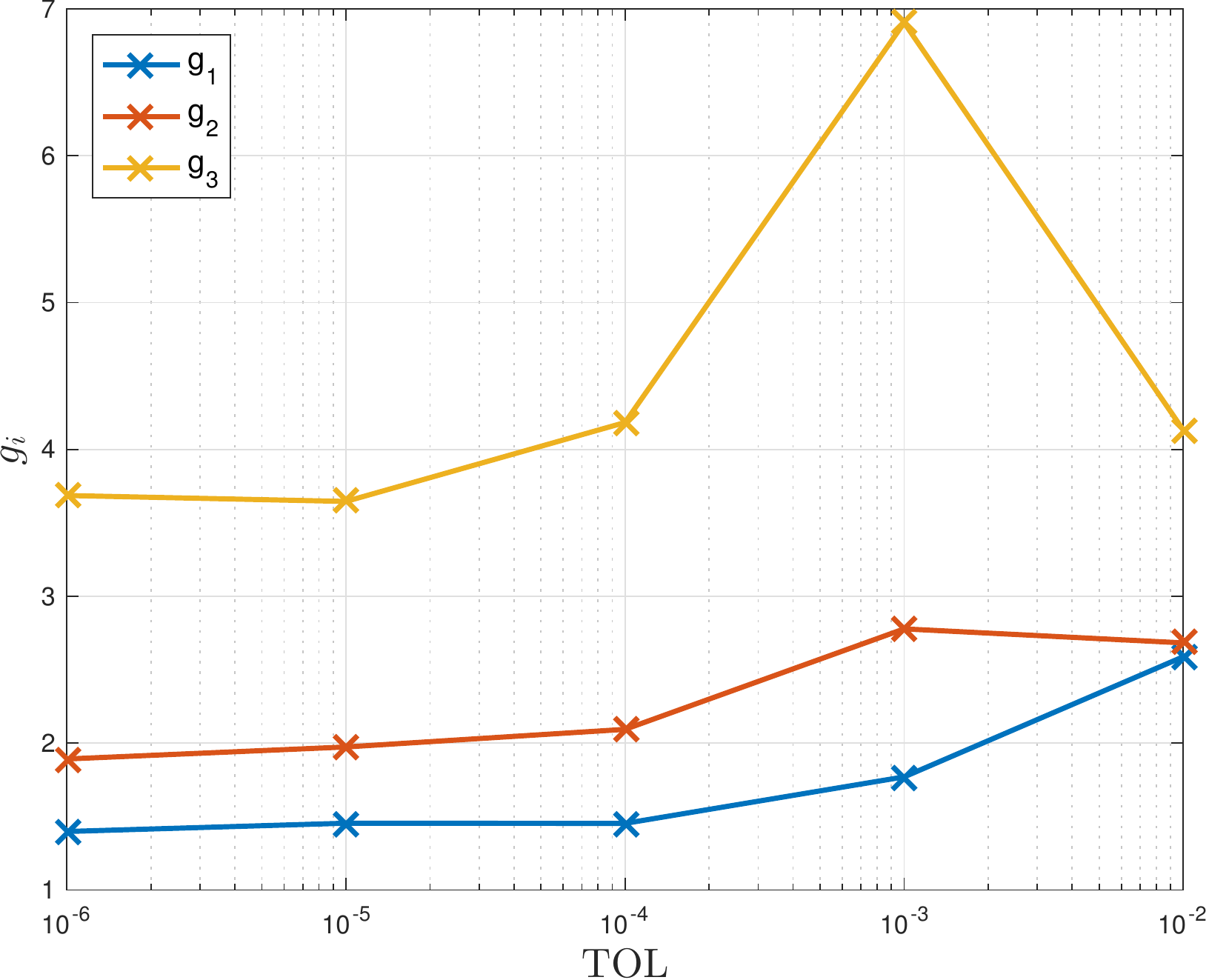}
   \caption{\La{Convergence of rates $g_i$ as computation proceeds.}}
   \label{fig:gi}
 \end{figure}

\subsection{Test 2 - $3d$ linear elasticity PDE with uncertain Lam\'e parameters}

\begin{figure}[t]
  \centering
  \includegraphics[width=0.27\linewidth]{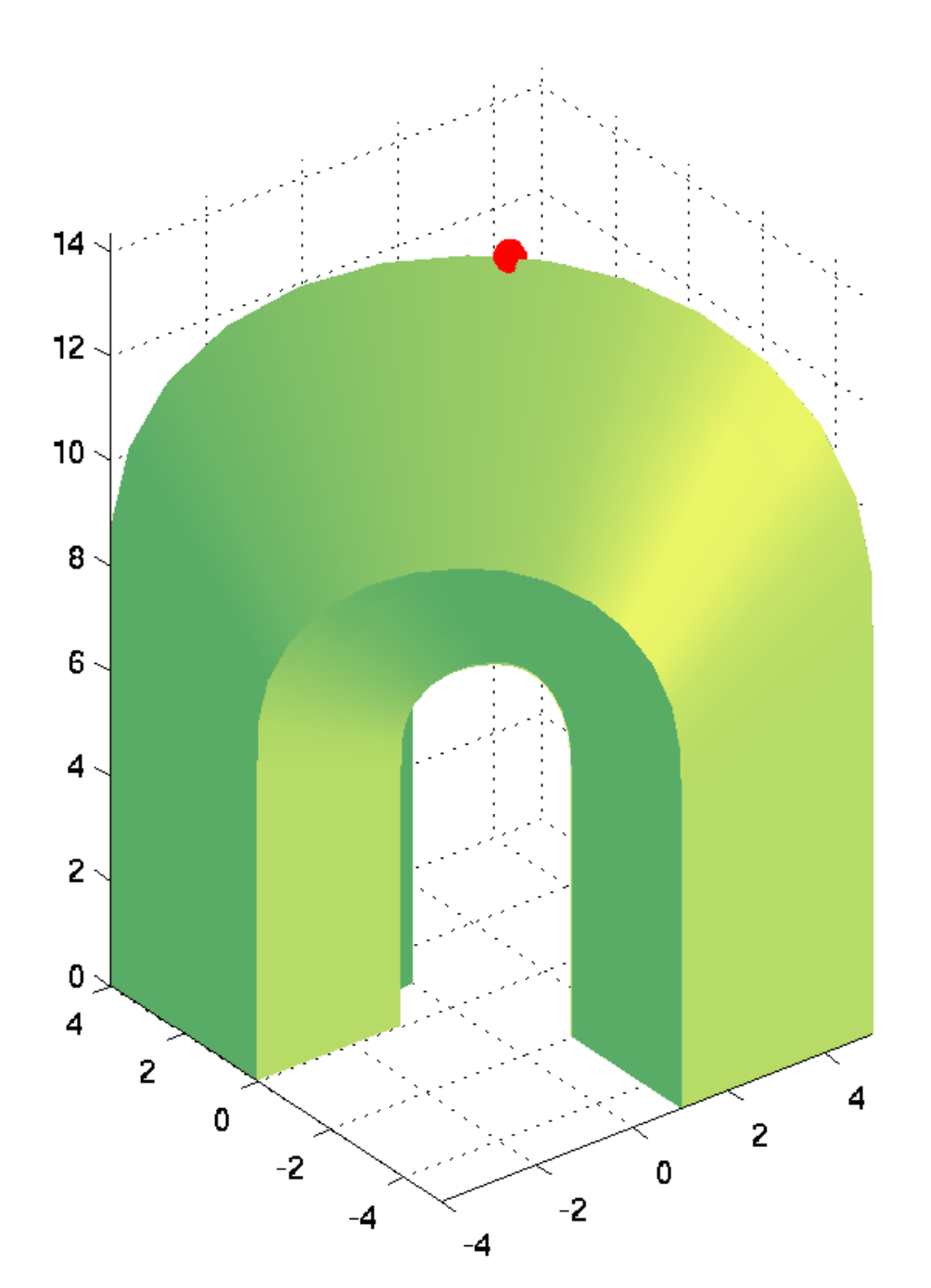}
  \includegraphics[width=0.27\linewidth]{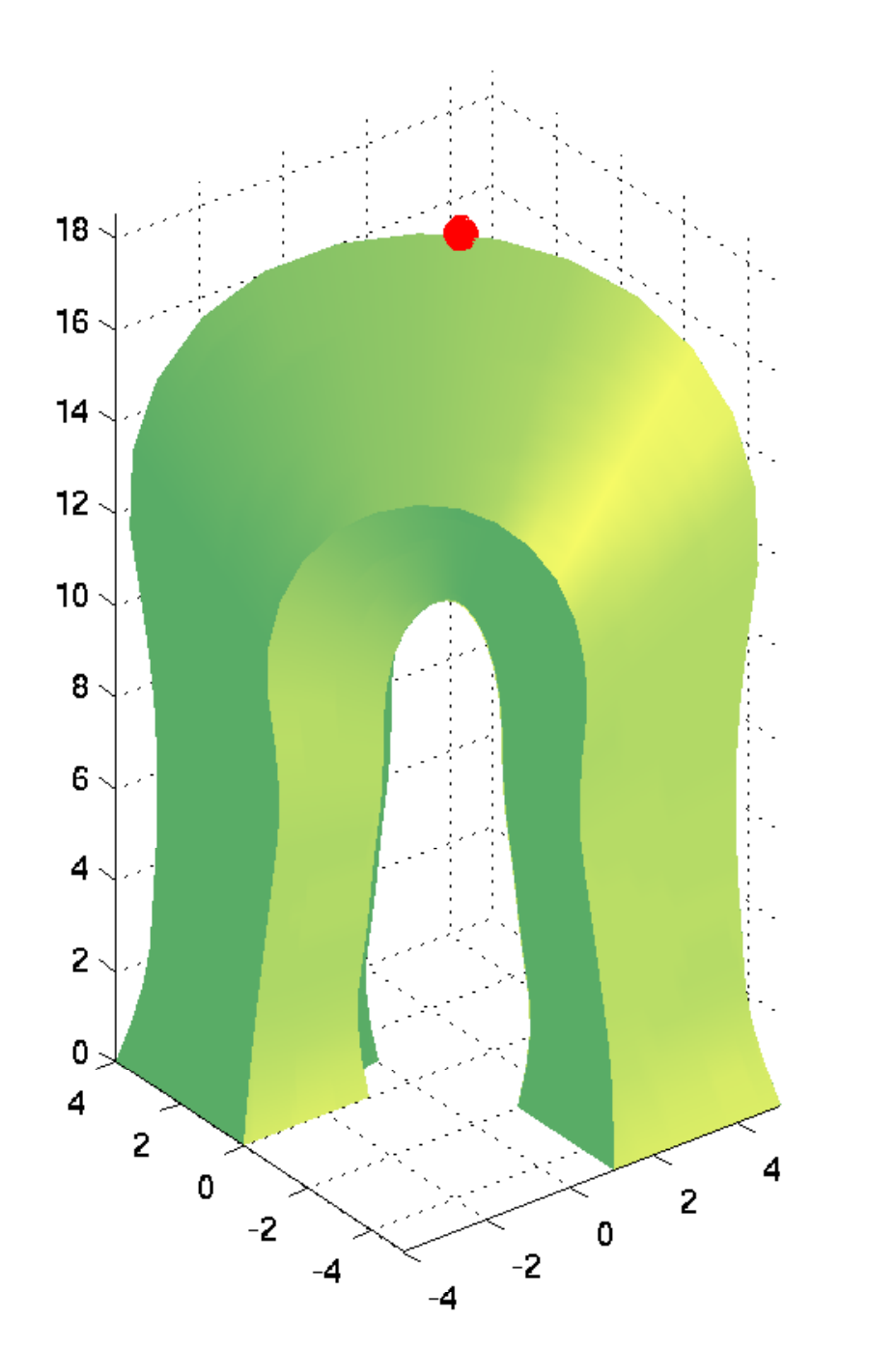}
  \caption{Test 2. Left: Undeformed geometry with $\mathbf{Q}$ marked by a red dot.
    Right: Deformation for a random value of the Lam\`e parameters.}\label{fig:test2_geometry_real}
\end{figure}

\begin{figure}[t]
  \centering
  \includegraphics[width=0.42\linewidth]{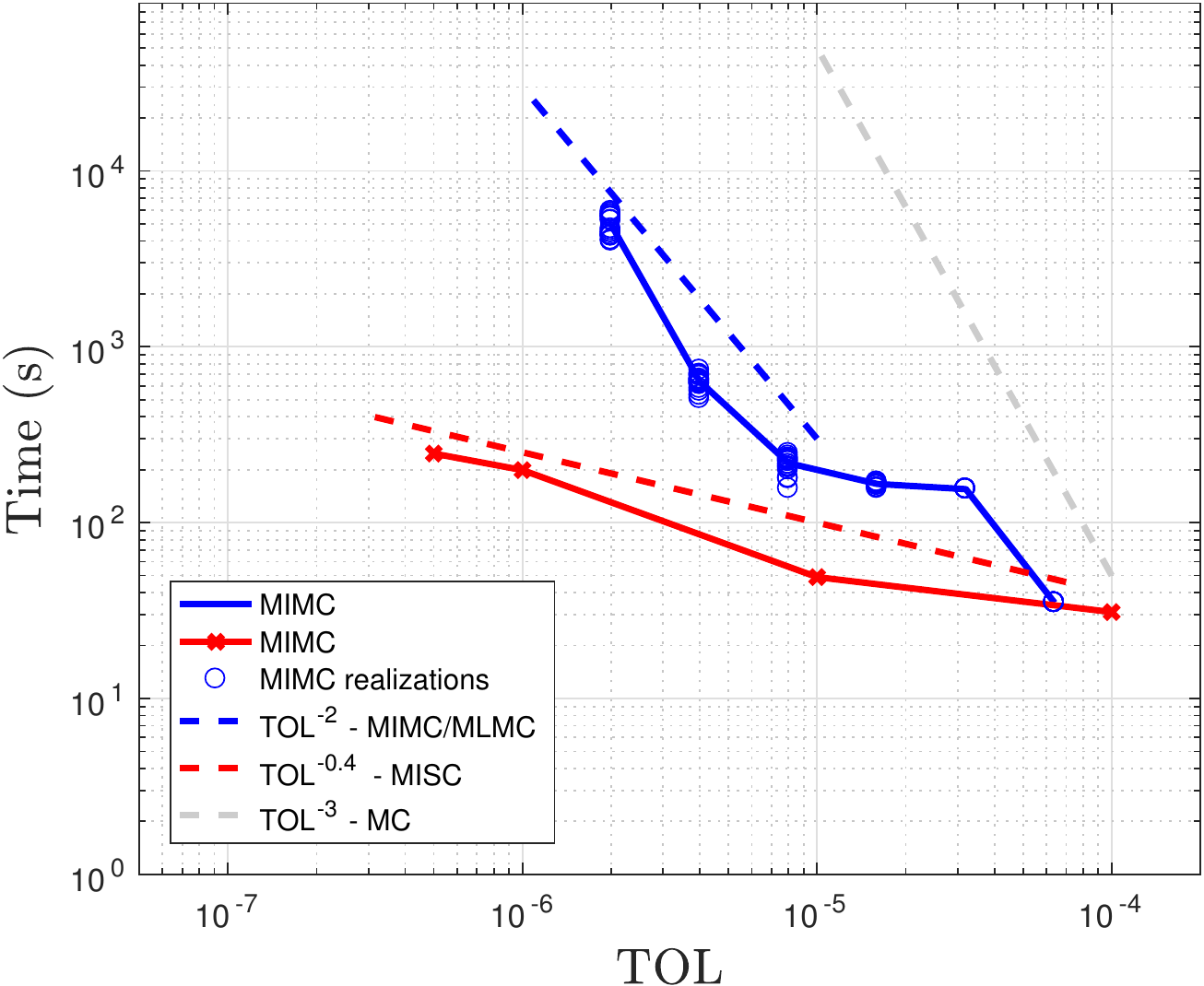}
  \includegraphics[width=0.42\linewidth]{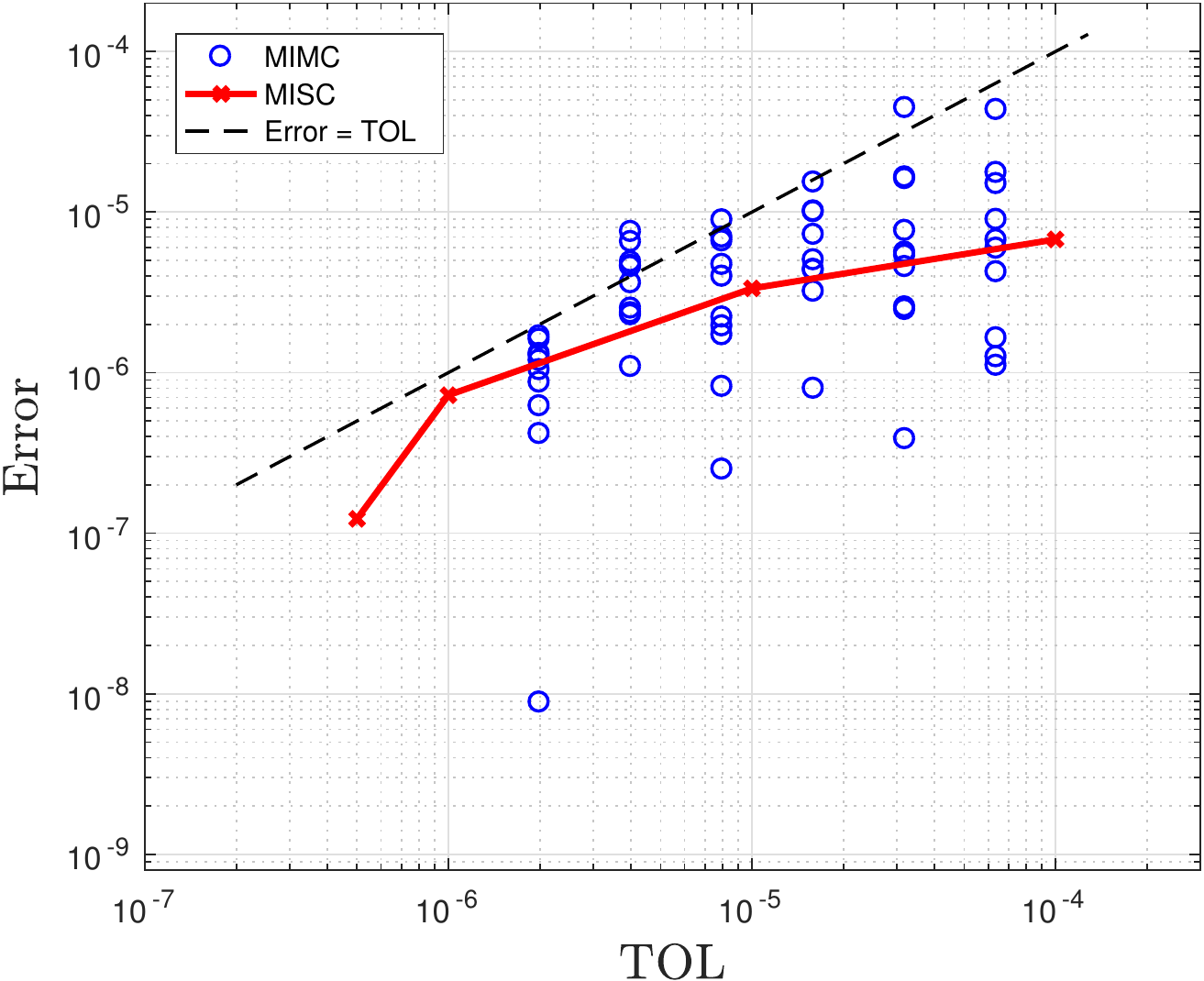}
  \caption{Test 2, computational results. Left: convergence results. Right: consistency of stopping criterion;
    the error has been computed against a sufficiently refined MISC solution.}\label{fig:test2_results}
\end{figure}

In this test we consider a slightly more complex problem than in the previous test, i.e., a linear elastic equation,
whose strong form reads
\[
  \begin{cases}
    -\text{div}[ \sigma( \uu(\xx,\yy) ) ] =\ff(\xx)  & \xx\in \mcB,\\
    \uu(\xx,\yy)=0 & \xx\in\partial \mcB_{\text{clamped}},\\
    \sigma( \uu(\xx,\yy) ) \cdot \nn =0 & \xx\in\partial \mcB_{\text{free}}.\\
  \end{cases}
\]
Although still elliptic in nature, this problem is computationally more demanding than the previous case because
the unknown is now the tri-dimensional displacement field $\uu: \mcB \rightarrow \Rset^3$.
The quantity $\sigma( \uu(\xx,\yy) )$ is the Cauchy stress tensor which, upon assuming
that the body is undergoing small deformations, can be related
to the displacement $\uu$ as
\[
\sigma(\uu(\xx,\yy)) = 2 \mu(\yy) \frac{\nabla \uu + [\nabla \uu]^T}{2} + \lambda(\yy) \text{div}(\uu) I,  
\]
where $\mu$, $\lambda$ are the Lam\'e constants and $I \in \Rset^{3 \times 3}$
is the identity matrix. The Lam\'e constants encode the mechanical properties of the material
and in this test they are assumed to be random variables, to model imperfect knowledge of such properties.
One notable practical example in which this might occur is 3d-printing, where the printer manufacturer guarantees
such properties only within a confidence range; see, e.g., \url{https://www.eos.info/material-m} .
In our experiments, we consider typical value ranges for titanium.
More specifically, we consider the following ranges for the Young's modulus $E$ and the Poission's ratio $\nu$, 
\[
  E \sim \mathcal{U}(105 \times 10^{9} \text{ Pa}, 120 \times 10^{9} \text{ Pa}), \quad \nu \sim \mathcal{U}(0.265,0.34),
\]
and then link these to the Lam\`e parameters by the well-known equations
\[
  \mu=\frac{E}{2(1+\nu)}, \quad \lambda=\frac{E \nu}{(1+\nu)(1-2\nu)}.
\]
Finally, we let $y_1,y_2 \sim \mathcal{U}([-1,1])$, and $y_1 \rightarrow E$, $y_2 \rightarrow \nu$
by linear maps. 
We consider the ``horse-shoe'' domain in Figure \ref{fig:test2_geometry_real}-left as computational domain $\mcB$.
The bottom end $z=0$ is kept fixed ($\mcB_{\text{clamped}}$ in the equation) and the rest of the body is free of any constraints
($\mcB_{\text{free}}$ in the equation). The body is pulled upward by a vertical force
$\ff=[0,\,0,\,10^{6}]$ N/m$^3$, and we are interested in computing the expected elongation
measured at $\mathbf{Q}$, the point marked by a red dot in Figure \ref{fig:test2_geometry_real}-left,
i.e.,  $\Phi(\uu) = \uu(\mathbf{Q})$;
an example of deformation obtained by two random values of the Lam\`e parameters
is shown in Figure \ref{fig:test2_geometry_real}-right
(magnified by a factor suitable to make it visible in a plot).
It is straightforward to see that the problem is well-posed for $\rho$-almost every $\yy \in \Gamma$ in the
vector-valued Hilbert space $V = [H^1(\mcB)]^3$.
As in Test 1, the IGA solver is set to use NURBS of degree $p=2$ with maximal continuity
everywhere except along the ridge of the horseshoe, and power-law-scaled knots in the parametric domain.
The rates $r_i, c_i$ are assessed numerically as 
$r_i=[2.5,2.5,2.5]$, $c_1=[1,1,1]$, implying convergence $\mathcal{O}(\text{TOL}^{-2})$
for MIMC and MLMC, $\mathcal{O}(\text{TOL}^{-3})$ for Monte Carlo
and $\mathcal{O}(\text{TOL}^{-0.4})$ up to preasymptotic and logarithmic terms for MISC.
Results are shown in Figure \ref{fig:test2_results},
and again confirm the predicted rates and show that MISC is significantly better than MIMC.
We omit the convergence results for MLMC, which, as in the previous test, is expected to converge with a
trend analogous to MIMC.

We conclude the discussion on this test by mentioning in-passing that in a sense
we are artificially increasing the complexity of the problem by using the Lam\`e
parameters instead of the Young's Modulus and Poisson's ratio.
Indeed, since both Lam\'e constants depend linearly on the Young Modulus $E$, 
the solution is inversely proportional to $E$, given the linearity of the PDE at hand.
Thus, $\phi(E,\nu) = 1/E \times \phi^*(\nu)$ and therefore
the MISC algorithm could be used over four indices (three in space and one for $\nu$) instead of five.
Nonetheless, we choose the formulation with five indices because our goal is to showcase
the computational efficiency of MISC in high-dimensional problems.
\La{However, we point out that even a MISC algorithm with adaptive selection of multi-index set
  would not recognize that there is ``one less inherent dimension'' in the problem,
  because the dependence is still non-linear in $E$.
  In other words, the fact that the true dependence is $\phi(E,\nu)=1/E \times \phi^*(\nu)$ implies that the adaptive algorithm
  would still need to refine with respect to $E$, and there is no particular structure of the multi-index set to be discovered:
  more specifically, the indices corresponding to mixed $E,\nu$ corrections are still required, because the size of these
  corrections are proportional to the mixed derivatives of $\phi$ with respect to $\nu$ and $E$, which are never zero.}

\section{Conclusions} \label{section:conclusions}

In this paper we extended the MISC methodology for solving elliptic PDEs with random coefficients
to non-square domains by using isogeometric analysis (IGA) solvers, which fit perfectly into the MISC framework due
to their tensor-structure construction.
\La{The MISC methodology aims at reaching a given accuracy in the approximation of the solution of the PDE
  at a minimal cost, by properly balancing the errors on the physical and stochastic domain.
  To do so, our implementation of MISC requires the knowledge of some parameters $r_i, c_i, g_i$,
  that are used to build the multi-index set prescribing the set of spatial and stochastic grids to be used.
  Note however that the multi-index set could also be chosen by a fully-adaptive a-posteriori scheme,
  which can be obtained e.g. with a straightforward extension of the dimension-adaptive
  scheme proposed in \cite{gerstner.griebel:adaptive,schillings.schwab:inverse,nobile.eal:adaptive-lognormal},
  which is classical in UQ literature.  
  We have proposed a practical algorithm with on-the-run tuning of parameters,
  and performed some tests which show the superiority of this approach
  with respect to the standard multi-level/multi-index Monte Carlo method.
  The superiority of MISC with respect to Monte-Carlo-based multi-level/multi-index methods is due to the fact that the
  stochastic sampling is now done with highly-efficient quadrature rules, which exploit the smoothness of the
  solution of the problem with respect to the random variables. Such smoothness is due to the properties of the
  differential operator (in this work, second-order elliptic), and some problems might not feature it (e.g., hyperbolic PDEs
  with shock formation); in such cases, Multi-Level/Multi-Index Monte Carlo should be used instead \cite{Pisaroni:MLMC-opt,mohamed:MOMC}.}
\La{The use of IGA solvers is important for MISC, since MISC needs a tensorized PDE solver to be fully effective;
  moreover, IGA solvers are advantageous in that they can represent exactly a wide range of geometries and they employ
  highly regular basis functions. However, IGA solvers are relatively young and still
  suffer some limitations when it comes to working on complex geometries. In any case,}
observe that in principle any tensor method able to deal with
non-square geometries in a tensorized fashion, such as Finite Differences, Finite Volumes, and $\mathbb{Q}^k$ finite element
on hexahedral meshes could be used as well. 
Research directions currently under investigations that could benefit from this approach are:
\begin{itemize}
\item forward UQ problems on domains with uncertain shape
  (and as a further step, shape optimization under uncertainty):
  indeed, the B-splines/NURBS representation of a geometry allows us to describe
  deviations from a nominal domain in a very straightforward manner.
\item UQ problems defined on unions of disjoint subdomains (``patches'' in the IGA literature).
  In this case, MISC could be allowed to choose different meshes in each
  subdomain, resulting in anisotropic meshes that would refine only a few of the subdomains.
  The resolution of the PDE on the corresponding non-conformal meshes could be performed by
  resorting to e.g. Lagrange multipliers \cite{Brivadis:mortar,Hesch:mortar,Temizer2012115} or discontinuous Galerkin methods
  \cite{antonietti:mortar,nitsche0,Arnold:DG} to enforce continuity at the interfaces of the subdomains. 
\end{itemize}

\section*{Acknowledgment}
The authors would like to thank the Isaac Newton Institute for
Mathematical Sciences, Cambridge, for support and hospitality during
the programme ``Uncertainty quantification for complex systems: theory
and methodologies'' supported by EPSRC grant no EP/K032208/1, where work on this paper was undertaken.
Part of this research was carried out while the authors visited the
Banff International Research Station for Mathematical Innovation and Discovery (BIRS),
for the workshop ``Computational Uncertainty Quantification'' in October 2017
(\url{https://www.birs.ca/events/2017/5-day-workshops/17w5072})
organized by Serge Prudhomme, Roger Ghanem, Mohammad Motamed, and Ra\'{u}l Tempone.
The hospitality and support of BIRS is acknowledged with gratitude.
This work was supported by the KAUST Office of Sponsored Research (OSR)
under award numbers URF/1/2281-01-01 and URF/1/2584-01-01 in the
KAUST Competitive Research Grants Program-Round 3 and 4, respectively.
Lorenzo Tamellini also received support from the European Union’s Horizon 2020
research and innovation program through the grant no. 680448 ``CAxMan'',
and by the GNCS 2018 project ``Metodi non conformi per equazioni alle derivate parziali''.

\section*{References}
\bibliographystyle{elsarticle-num}
\bibliography{IGA_biblio,UQ_biblio,sparse_grids_biblio,multipatch_biblio}

\def\cprime{$'$}
\begin{thebibliography}{10}
\expandafter\ifx\csname url\endcsname\relax
  \def\url#1{\texttt{#1}}\fi
\expandafter\ifx\csname urlprefix\endcsname\relax\def\urlprefix{URL }\fi
\expandafter\ifx\csname href\endcsname\relax
  \def\href#1#2{#2} \def\path#1{#1}\fi

\bibitem{ghanem:UQbook}
R.~Ghanem, D.~Higdon, H.~Owhadi, Handbook of Uncertainty Quantification,
  Handbook of Uncertainty Quantification, Springer International Publishing,
  2016.

\bibitem{sullivan:UQbook}
T.~Sullivan, Introduction to Uncertainty Quantification, Texts in Applied
  Mathematics, Springer International Publishing, 2015.

\bibitem{smith:UQbook}
R.~Smith, Uncertainty Quantification: Theory, Implementation, and Applications,
  Computational Science and Engineering, Society for Industrial and Applied
  Mathematics, 2013.

\bibitem{nobile.tempone.eal:aniso}
F.~Nobile, R.~Tempone, C.~Webster, An anisotropic sparse grid stochastic
  collocation method for partial differential equations with random input data,
  SIAM J. Numer. Anal. 46~(5) (2008) 2411--2442.

\bibitem{nobile.eal:optimal-sparse-grids}
F.~Nobile, L.~Tamellini, R.~Tempone, Convergence of quasi-optimal sparse-grid
  approximation of {H}ilbert-space-valued functions: application to random
  elliptic {PDE}s, Numerische {M}athematik 134~(2) (2016) 343--388.

\bibitem{cohen.devore.schwab:nterm2}
A.~Cohen, R.~Devore, C.~Schwab, Analytic regularity and polynomial
  approximation of parametric and stochastic elliptic {PDE'S}, Anal. Appl.
  (Singap.) 9~(1) (2011) 11--47.

\bibitem{chkifa:adaptive-taylor}
A.~Chkifa, A.~Cohen, R.~Devore, C.~Schwab, Sparse adaptive {T}aylor
  approximation algorithms for parametric and stochastic elliptic {PDE}s,
  ESAIM: Mathematical Modelling and Numerical Analysis 47~(1) (2013) 253--280.

\bibitem{scheichl.giles:MLMC}
K.~Cliffe, M.~Giles, R.~Scheichl, A.~Teckentrup, Multilevel monte carlo methods
  and applications to elliptic pdes with random coefficients, Computing and
  Visualization in Science 14~(1) (2011) 3--15.

\bibitem{hajiali.eal:MultiIndexMC}
A.-L. Haji-Ali, F.~Nobile, R.~Tempone, {M}ulti-{i}ndex {M}onte {C}arlo: when
  sparsity meets sampling, Numerische Mathematik (2015) 1--40.

\bibitem{peherstorfer:MFsurvey}
B.~Peherstorfer, K.~Willcox, M.~Gunzburger, Survey of multifidelity methods in
  uncertainty propagation, inference, and optimization, SIAM Review 60~(3)
  (2018) 550--591.

\bibitem{hajiali.eal:MISC1}
A.~Haji-Ali, F.~Nobile, L.~Tamellini, R.~Tempone, Multi-index stochastic
  collocation for random \{PDEs\}, Computer Methods in Applied Mechanics and
  Engineering 306 (2016) 95 -- 122.

\bibitem{hajiali.eal:MISC2}
A.-L. Haji-Ali, F.~Nobile, L.~Tamellini, R.~Tempone, Multi-index {S}tochastic
  {C}ollocation convergence rates for random {PDE}s with parametric regularity,
  Foundations of {C}omputational {M}athematics 16~(6) (2016) 1555--1605.

\bibitem{Bungartz.Griebel.Roschke.ea:pointwise.conv}
H.-J. Bungartz, M.~Griebel, D.~R\"oschke, C.~Zenger, Pointwise convergence of
  the combination technique for the {L}aplace equation, East-West J. Numer.
  Math. 2 (1994) 21--45.

\bibitem{b.griebel:acta}
H.~Bungartz, M.~Griebel, Sparse grids, Acta Numer. 13 (2004) 147--269.

\bibitem{Griebel.schneider.zenger:combination}
M.~Griebel, M.~Schneider, C.~Zenger, A combination technique for the solution
  of sparse grid problems, in: P.~de~Groen, R.~Beauwens (Eds.), {Iterative
  Methods in Linear Algebra}, IMACS, Elsevier, North Holland, 1992, pp.
  263--281.

\bibitem{Hegland:combination}
M.~Hegland, J.~Garcke, V.~Challis, The combination technique and some
  generalisations, Linear Algebra and its Applications 420~(2–3) (2007) 249
  -- 275.

\bibitem{smolyak:quadrature}
S.~Smolyak, Quadrature and interpolation formulas for tensor products of
  certain classes of functions, Dokl. Akad. Nauk SSSR 4 (1963) 240--243.

\bibitem{barthelmann.novak.ritter:high}
V.~Barthelmann, E.~Novak, K.~Ritter, High dimensional polynomial interpolation
  on sparse grids, Adv. Comput. Math. 12~(4) (2000) 273--288.

\bibitem{quarteroni.sacco.eal:numerical}
A.~Quarteroni, R.~Sacco, F.~Saleri, Numerical mathematics, 2nd Edition, Vol.~37
  of Texts in Applied Mathematics, Springer-Verlag, Berlin, 2007.

\bibitem{teckentrup.etal:MLSC}
A.~L. Teckentrup, P.~Jantsch, C.~G. Webster, M.~Gunzburger, A {M}ultilevel
  {S}tochastic {C}ollocation {M}ethod for {P}artial {D}ifferential {E}quations
  with {R}andom {I}nput {D}ata, SIAM/ASA Journal on Uncertainty Quantification
  3~(1) (2015) 1046--1074.

\bibitem{van-wyk:MLSC}
H.~W. van Wyk, Multilevel sparse grid methods for elliptic partial differential
  equations with random coefficients, arXiv arXiv:1404.0963, e-print (2014).

\bibitem{hps13}
H.~Harbrecht, M.~Peters, M.~Siebenmorgen, On multilevel quadrature for elliptic
  stochastic partial differential equations, in: Sparse Grids and Applications,
  Vol.~88 of Lecture Notes in Computational Science and Engineering, Springer,
  2013, pp. 161--179.

\bibitem{babuska.nobile.eal:stochastic2}
I.~Babu\v{s}ka, F.~Nobile, R.~Tempone, A stochastic collocation method for
  elliptic partial differential equations with random input data, SIAM Review
  52~(2) (2010) 317--355.

\bibitem{xiu.hesthaven:high}
D.~Xiu, J.~Hesthaven, High-order collocation methods for differential equations
  with random inputs, SIAM J. Sci. Comput. 27~(3) (2005) 1118--1139.

\bibitem{Hughes:2005}
T.~Hughes, J.~Cottrell, Y.~Bazilevs, Isogeometric analysis: {CAD}, finite
  elements, {NURBS}, exact geometry and mesh refinement, Computer Methods in
  Applied Mechanics and Engineering 194~(39) (2005) 4135--4195.

\bibitem{IGA-book}
J.~A. Cottrell, T.~J.~R. Hughes, Y.~Bazilevs, Isogeometric {A}nalysis: toward
  integration of {CAD} and {FEA}, John Wiley \& Sons, 2009.

\bibitem{acta-IGA}
L.~Beirao Da~Veiga, A.~Buffa, G.~Sangalli, R.~V\'{a}zquez, Mathematical
  analysis of variational isogeometric methods, Acta Numerica 23 (2014)
  157--287.

\bibitem{Benzaken20171215}
J.~Benzaken, A.~Herrema, M.-C. Hsu, J.~Evans, A rapid and efficient
  isogeometric design space exploration framework with application to
  structural mechanics, Computer Methods in Applied Mechanics and Engineering
  316 (2017) 1215 -- 1256, special Issue on Isogeometric Analysis: Progress and
  Challenges.

\bibitem{manzoni.heltai:RB-IGA}
A.~Manzoni, F.~Salmoiraghi, L.~Heltai, Reduced basis isogeometric methods
  (rb-iga) for the real-time simulation of potential flows about parametrized
  \{NACA\} airfoils, Computer Methods in Applied Mechanics and Engineering 284
  (2015) 1147 -- 1180, isogeometric Analysis Special Issue.

\bibitem{WILHELM2016}
M.~Wilhelm, L.~Ded\`{e}, L.~M. Sangalli, P.~Wilhelm, Igs: An isogeometric
  approach for smoothing on surfaces, Computer Methods in Applied Mechanics and
  Engineering 302 (2016) 70 -- 89.

\bibitem{corno_UQIGA}
J.~{Corno}, C.~{de Falco}, H.~{De Gersem}, S.~{Schöps}, Isogeometric analysis
  simulation of tesla cavities under uncertainty, in: 2015 International
  Conference on Electromagnetics in Advanced Applications (ICEAA), 2015, pp.
  1508--1511.

\bibitem{LI_UQIGA}
K.~Li, W.~Gao, D.~Wu, C.~Song, T.~Chen, Spectral stochastic isogeometric
  analysis of linear elasticity, Computer Methods in Applied Mechanics and
  Engineering 332 (2018) 157 -- 190.

\bibitem{HIEN2017}
T.~D. Hien, H.-C. Noh, Stochastic isogeometric analysis of free vibration of
  functionally graded plates considering material randomness, Computer Methods
  in Applied Mechanics and Engineering 318 (2017) 845 -- 863.

\bibitem{RAHMAN2018}
S.~Rahman, A galerkin isogeometric method for karhunen–loève approximation
  of random fields, Computer Methods in Applied Mechanics and Engineering 338
  (2018) 533 -- 561.

\bibitem{beck.eal:sparse-IGA}
J.~Beck, G.~Sangalli, L.~Tamellini, A sparse-grid isogeometric solver,
  {Computer Methods in Applied Mechanics and Engineering} 335~(--) (2018)
  128--151.

\bibitem{farin2001curves}
G.~Farin, Curves and Surfaces for CAGD: A Practical Guide, The Morgan Kaufmann
  Series in Computer Graphics, Elsevier Science, 2001.

\bibitem{cohen2001geometric}
E.~Cohen, R.~Riesenfeld, G.~Elber, Geometric modeling with splines: an
  introduction, Vol.~1, AK Peters Wellesley, MA, 2001.

\bibitem{montardini:collocation}
M.~Montardini, G.~Sangalli, L.~Tamellini, Optimal-order isogeometric
  collocation at {G}alerkin superconvergent points, Computer Methods in Applied
  Mechanics and Engineering.

\bibitem{gomez2016variational}
H.~Gomez, L.~De~Lorenzis, The variational collocation method, Computer Methods
  in Applied Mechanics and Engineering 309 (2016) 152--181.

\bibitem{anitescu2015isogeometric}
C.~Anitescu, Y.~Jia, Y.~J. Zhang, T.~Rabczuk, An isogeometric collocation
  method using superconvergent points, Computer Methods in Applied Mechanics
  and Engineering 284 (2015) 1073--1097.

\bibitem{casquero2016isogeometric}
H.~Casquero, L.~Liu, Y.~Zhang, A.~Reali, H.~Gomez, Isogeometric collocation
  using analysis-suitable t-splines of arbitrary degree, Computer Methods in
  Applied Mechanics and Engineering 301 (2016) 164--186.

\bibitem{DOKKEN2018}
T.~Dokken, V.~Skytt, O.~Barrowclough, Trivariate spline representations for
  computer aided design and additive manufacturing, Computers \& Mathematics
  with Applications.

\bibitem{MASSARWI2016}
F.~Massarwi, G.~Elber, A b-spline based framework for volumetric object
  modeling, Computer-Aided Design 78 (2016) 36 -- 47, sPM 2016.

\bibitem{KUDELA2016406}
L.~Kudela, N.~Zander, S.~Kollmannsberger, E.~Rank, Smart octrees: Accurately
  integrating discontinuous functions in 3d, Computer Methods in Applied
  Mechanics and Engineering 306 (2016) 406 -- 426.

\bibitem{Marussig2018}
B.~Marussig, T.~J.~R. Hughes, A review of trimming in isogeometric analysis:
  Challenges, data exchange and simulation aspects, Archives of Computational
  Methods in Engineering 25~(4) (2018) 1059--1127.

\bibitem{ernst.eal:collocation-logn}
O.~G. Ernst, B.~Sprungk, L.~Tamellini, {Convergence of Sparse Collocation for
  Functions of Countably Many Gaussian Random Variables (with Application to
  Lognormal Elliptic Diffusion Problems)}, {SIAM Journal on Numerical Analysis}
  56~(2) (2018) 877--905.

\bibitem{trefethen:comparison}
L.~N. Trefethen, Is {G}auss quadrature better than {C}lenshaw-{C}urtis?, SIAM
  Rev. 50~(1) (2008) 67--87.

\bibitem{nobile.etal:leja}
F.~Nobile, L.~Tamellini, R.~Tempone, Comparison of {C}lenshaw–{C}urtis and
  {L}eja {Q}uasi-{O}ptimal {S}parse {G}rids for the {A}pproximation of {R}andom
  {PDE}s, in: R.~M. Kirby, M.~Berzins, J.~S. Hesthaven (Eds.), Spectral and
  High Order Methods for Partial Differential Equations - ICOSAHOM '14, Vol.
  106 of Lecture Notes in Computational Science and Engineering, Springer
  International Publishing, 2015, pp. 475--482.

\bibitem{narayan:Leja}
A.~Narayan, J.~D. Jakeman, Adaptive {L}eja {S}parse {G}rid {C}onstructions for
  {S}tochastic {C}ollocation and {H}igh-{D}imensional {A}pproximation, SIAM
  Journal on Scientific Computing 36~(6) (2014) A2952--A2983.

\bibitem{Chkifa:leja}
A.~Chkifa, On the {L}ebesgue constant of {L}eja sequences for the complex unit
  disk and of their real projection, Journal of Approximation Theory 166~(0)
  (2013) 176 -- 200.

\bibitem{back.nobile.eal:comparison}
J.~B\"ack, F.~Nobile, L.~Tamellini, R.~Tempone, Stochastic spectral {G}alerkin
  and collocation methods for {PDE}s with random coefficients: a numerical
  comparison, in: Spectral and High Order Methods for Partial Differential
  Equations, Vol.~76 of Lecture Notes in Computational Science and Engineering,
  Springer, 2011, pp. 43--62.

\bibitem{bieri:sparse.tensor.coll}
M.~Bieri, A sparse composite collocation finite element method for elliptic
  spdes., SIAM Journal on Numerical Analysis 49~(6) (2011) 2277--2301.

\bibitem{gerstner.griebel:adaptive}
T.~Gerstner, M.~Griebel, Dimension-adaptive tensor-product quadrature,
  Computing 71~(1) (2003) 65--87.

\bibitem{schillings.schwab:inverse}
C.~Schillings, C.~Schwab, Sparse, adaptive {S}molyak quadratures for {B}ayesian
  inverse problems, Inverse Problems 29~(6).

\bibitem{nobile.eal:adaptive-lognormal}
F.~Nobile, L.~Tamellini, F.~Tesei, R.~Tempone, An adaptive sparse grid
  algorithm for elliptic {PDE}s with lognormal diffusion coefficient, in:
  J.~Garcke, D.~Pfl\"uger (Eds.), Sparse Grids and Applications -- Stuttgart
  2014, Vol. 109 of Lecture Notes in Computational Science and Engineering,
  Springer International Publishing Switzerland, 2016, pp. 191--220.

\bibitem{Guo1986}
B.~Guo, I.~Babu{\v{s}}ka, The h-p version of the finite element method,
  Computational Mechanics 1~(3) (1986) 203--220.

\bibitem{Babuska_book}
I.~Babu{\v{s}}ka, T.~Strouboulis, The finite element method and its
  reliability, Numerical Mathematics and Scientific Computation, The Clarendon
  Press Oxford University Press, New York, 2001.

\bibitem{SERSH13}
D.~Schillinger, J.~Evans, A.~Reali, M.~Scott, T.~Hughes, Isogeometric
  collocation: cost comparison with {G}alerkin methods and extension to
  adaptive hierarchical {NURBS} discretizations, Comput. Methods Appl. Mech.
  Engrg. 267 (2013) 170 -- 232.

\bibitem{Guignard:a-post}
D.~S. Guignard, F.~Nobile, A posteriori error estimation for the stochastic
  collocation finite element method, Mathicse Report nr 24.2017.

\bibitem{mohamed:MOMC}
M.~Motamed, D.~Appel{\"o}, A multi order discontinuous galerkin monte carlo
  method for hyperbolic problems with stochastic parameters, SIAM Journal on
  Numerical Analysis 56~(1) (2018) 448--468.

\bibitem{VAZQUEZ2016523}
R.~Vazquez, A new design for the implementation of isogeometric analysis in
  octave and matlab: Geopdes 3.0, Computers \& Mathematics with Applications
  72~(3) (2016) 523 -- 554.

\bibitem{cohen_devore_2015}
A.~Cohen, R.~DeVore, Approximation of high-dimensional parametric pdes, Acta
  Numerica 24 (2015) 1–159.

\bibitem{back.nobile.eal:lognormal}
J.~Beck, F.~Nobile, L.~Tamellini, R.~Tempone, A {Q}uasi-optimal {S}parse
  {G}rids {P}rocedure for {G}roundwater {F}lows, in: Spectral and High Order
  Methods for Partial Differential Equations - ICOSAHOM 2012, Vol.~95 of
  Lecture Notes in Computational Science and Engineering, Springer, 2014, pp.
  1--16.

\bibitem{ghanem.spanos:book}
R.~G. Ghanem, P.~D. Spanos, Stochastic finite elements: a spectral approach,
  Springer-Verlag, New York, 1991.

\bibitem{chkifa:nonlinear}
A.~Chkifa, A.~Cohen, C.~Schwab, Breaking the curse of dimensionality in sparse
  polynomial approximation of parametric {PDEs}, Journal de Mathématiques
  Pures et Appliquées 103~(2) (2015) 400 -- 428.

\bibitem{hajiali.eal:continuationMLMC}
N.~Collier, A.~Haji-Ali, F.~Nobile, E.~von Schwerin, R.~Tempone, A continuation
  multilevel monte carlo algorithm, BIT Numerical Mathematics (2014) 1--34.

\bibitem{Pisaroni:MLMC-opt}
M.~Pisaroni, F.~Nobile, P.~Leyland, A multilevel monte carlo evolutionary
  algorithm for robust aerodynamic shape design, in: 18th AIAA/ISSMO
  Multidisciplinary Analysis and Optimization Conference. Denver, Colorado,
  2017.

\bibitem{Brivadis:mortar}
E.~Brivadis, A.~Buffa, B.~Wohlmuth, L.~Wunderlich, Isogeometric mortar methods,
  Computer Methods in Applied Mechanics and Engineering 284 (2015) 292 -- 319,
  isogeometric Analysis Special Issue.

\bibitem{Hesch:mortar}
C.~Hesch, P.~Betsch, Isogeometric analysis and domain decomposition methods,
  Computer Methods in Applied Mechanics and Engineering 213–216 (2012) 104 --
  112.

\bibitem{Temizer2012115}
{\.I}.~Temizer, P.~Wriggers, T.~Hughes, Three-dimensional mortar-based
  frictional contact treatment in isogeometric analysis with {NURBS}, Comput.
  Methods Appl. Mech. Engrg. 209--212 (2012) 115 -- 128.

\bibitem{antonietti:mortar}
P.~Antonietti, I.~Mazzieri, A.~Quarteroni, F.~Rapetti, Non-conforming high
  order approximations of the elastodynamics equation, Computer Methods in
  Applied Mechanics and Engineering 209–212 (2012) 212 -- 238.

\bibitem{nitsche0}
J.~Nitsche, \"{U}ber ein {V}ariationsprinzip zur {L}\"osung von
  {D}irichlet-{P}roblemen bei {V}erwendung von {T}eilr\"aumen, die keinen
  {R}andbedingungen unterworfen sind, Abh. Math. Sem. Univ. Hamburg 36 (1971)
  9--15, collection of articles dedicated to Lothar Collatz on his sixtieth
  birthday.

\bibitem{Arnold:DG}
D.~N. Arnold, F.~Brezzi, B.~Cockburn, L.~D. Marini, Unified analysis of
  discontinuous galerkin methods for elliptic problems, SIAM Journal on
  Numerical Analysis 39~(5) (2002) 1749--1779.

\end{thebibliography}

\end{document}